# RATES FOR BRANCHING PARTICLE APPROXIMATIONS OF CONTINUOUS-DISCRETE FILTERS

By Michael A. Kouritzin and Wei Sun

*University of Alberta and Concordia University*

Herein, we analyze an efficient branching particle method for asymptotic solutions to a class of continuous-discrete filtering problems. Suppose that $t \to X_t$ is a Markov process and we wish to calculate the measure-valued process $t \to \mu_t(\cdot) \doteq P\{X_t \in \cdot | \sigma\{Y_{t_k}, t_k \le t\}\}$, where $t_k = k\varepsilon$ and $Y_{t_k}$ is a distorted, corrupted, partial *observation* of $X_{t_k}$. Then, one constructs a particle system with observation-dependent branching and $n$ initial particles whose empirical measure at time $t$, $\mu_t^n$, closely approximates $\mu_t$. Each particle evolves independently of the other particles according to the law of the *signal* between observation times $t_k$, and branches with *small probability* at an observation time. For filtering problems where $\varepsilon$ is very small, using the algorithm considered in this paper requires far fewer computations than other algorithms that branch or interact all particles regardless of the value of $\varepsilon$. We analyze the algorithm on Lévy-stable signals and give rates of convergence for $E^{1/2}\{\|\mu_t^n - \mu_t\|_\gamma^2\}$, where $\|\cdot\|_\gamma$ is a Sobolev norm, as well as related convergence results.

**1. Introduction.** The filtering problems in many key, contemporary fields such as mathematical finance and communication networks initially appear to be resolved by the celebrated mathematical solutions of the Duncan–Mortensen–Zakai and Kushner–Stratonovich equations, which have been known for over three decades. However, upon further reflection, one realizes that these equations are neither computer workable nor applicable at large. More theory is required keeping: (a) the ultimate computer enduse, and (b) some real world applications in mind. Many of the corresponding

Received August 2003; revised April 2005.

[1]Supported by NSERC, PIMS, Lockheed Martin Naval Electronics and Surveillance Systems, Lockheed Martin Canada and Acoustic Positioning Research Inc., through a MITACS centre of excellence entitled “Prediction in Interacting Systems.”

filtering problems are large enough that the mere storage of the exact solutions is impractical. We require more implementable, practical methods of filtering, where the solutions are *almost* optimal and can be stored. The introduction of particle approximations is natural under these criteria.

The general problem of continuous-discrete filtering for Markov processes is concerned with extracting information about a continuous-time Markov process $t \to X_t$ called the *signal* based on the current record of discrete-time *observations* $\{Y_{t_k}, t_k \leq t\}$ that are probabilistically linked to the signal. The goal of filtering is to estimate past, present or future values of $\varphi(X_t)$ based on our observation record $\{Y_{t_k}, t_k \leq t\}$. Direct implementation of the mathematical solution to these filtering problems usually requires the on-line solution of an infinite-dimensional (often parabolic) equation (see, however, [11] for counter examples where such infinite-dimensional equation solution is not required), which is impossible to either implement precisely or store. For these reasons, one may be forced to approximate. One exciting method of approximation for continuous-discrete filtering problems was recently studied by Del Moral and collaborators (see [7] for one of the earlier works), where, instead of solving a parabolic equation on-line, one simulates particles so that the empirical measure of the particles is a good approximation to the solution of the differential equation. Then, to account for the incoming observations, one allows the particles to redistribute themselves to locations favored by the observations. This second branching or interacting step is devised to ensure that new information obtained through the observations can be incorporated into our conditional probability law of the signal given the observation record. A thorough account of this interesting interacting particle method can be found in [9].

More recently, algorithms have been considered in [1] and [8] that do not disturb most particles at each observation time and thereby introduce far less resampling noise. Indeed, the huge performance gained by only resampling those particles that need to be resampled was quantified experimentally in the former paper and theoretically in the latter. Herein, we further develop and study the cautious branching particle approach in [1], which was motivated in part by the particle system approximation scheme suggested by Sherman and Peskin [17] for the deterministic reaction-diffusion equations and by the earlier branching particle method of Crisan and Lyons [6]. To make our presentation clear, we choose to introduce and analyze our method on Lévy-stable signal processes, however, this particle approximation method is extendable well beyond our current setting as experiments have demonstrated.

Lévy-stable processes are one of the most basic and important classes of Markov processes. They are widely used in various economic and physical systems. In particular, the use of Lévy-stable processes in mathematical finance and communication networks has recently become more popular.

For instance, Lévy-stable models have been applied in the fields of portfolio theory, asset and option pricing (cf. [3, 14] and the references therein); and Lévy-stable processes have been used to model teletraffic and to approximate network traffic (cf. [10, 15]). These vital applications are motivation for us to analyze our method on Lévy-stable signals.

We let $(\Omega, \mathcal{F}, P^0)$ be a complete probability space and $E^0$ be the expectation with respect to $P^0$. Suppose that $X$ is a $\mathbb{R}^{d_1}$-valued Lévy-stable process on $(\Omega, \mathcal{F}, P^0)$ with index $\alpha \in (0,2]$ and spectral measure $\Gamma$ (cf. [16]), that is, $X$ is a stochastic process on $\mathbb{R}^{d_1}$ such that $X$ has independent increments, and there exists a finite measure $\Gamma$ on the unit sphere $S_{d_1}$ of $\mathbb{R}^{d_1}$ such that, for any $\theta = (\theta_1, \dots, \theta_{d_1})' \in \mathbb{R}^{d_1}$ and $0 \le s < t < \infty$,

$$\ln E^0\{\exp(i\theta'(X_t - X_s))\}$$
$$= \begin{cases} -(t-s)\displaystyle\int_{S_{d_1}} |\theta' z|^\alpha \left(1 - i\,\mathrm{sign}(\theta' z)\tan\left(\frac{\alpha\pi}{2}\right)\right)\Gamma(dz), & \text{for } \alpha \neq 1, \\ -(t-s)\displaystyle\int_{S_{d_1}} |\theta' z|\left(1 + \frac{2i}{\pi}\mathrm{sign}(\theta' z)\ln|\theta' z|\right)\Gamma(dz), & \text{for } \alpha = 1. \end{cases}$$

Hereafter, we use " $'$ " to denote the transpose of a vector. We let $0 < \varepsilon \le 1$, define $t_k \doteq k\varepsilon$ for $k = 1, 2, \dots$, and suppose that $V$ is a standard $\mathbb{R}^{d_2}$-valued Brownian motion on $(\Omega, \mathcal{F}, P^0)$ independent of $X$. Then, we consider calculating the conditional probability law of signal $X_t$ given the multi-dimensional observations $\{Y_{t_k}, t_k \le t\}$, defined by

$$Y_{t_k} = Y_{t_{k-1}} + h(X_{t_k})(t_k - t_{k-1}) + (V_{t_k} - V_{t_{k-1}}),$$

via change of measure and *particle approximation.*

Our particle approximation scheme can be summarized as follows: We consider a branching particle system which starts off with $n$ particles and each particle has the "opportunity" to branch and die every $\varepsilon$ seconds. A particle reaching $x$ at time $t_k-$ branches with small probability and in this unlikely event that the particle does branch, it either just dies or is replaced by two or more independent particles starting at $(t_k, x)$. Efficiencies are gained at observation times in two ways: The vast majority of particles do not branch at branching times for small $\varepsilon$, which reduces computation related to duplicating or removing particles, and branching decisions only depend on the very particle that may or may not branch so decisions require little processing. On the other hand, the number of particles in our scheme does not stay constant, but rather is a nontrivial martingale. Still, there are effective ways to control the number of particles in practice, by introducing additional births and deaths that do not bias estimates, and thereby to keep the computations essentially constant over the various observation times.

Suppose that $\delta_x$ denotes the Dirac delta measure at $x$ and

$$\mu_t^n \doteq \frac{1}{n} \sum_{i=1}^{\|\mu_t^n\|} \delta_{X_t^{i,n}} \tag{1}$$

is the empirical measure of the particle system if there are $\|\mu_t^n\|$ particles $\{X_t^{1,n}, \dots, X_t^{\|\mu_t^n\|,n}\}$ alive at time $t$. Then, among other things, our results will imply that

$$\left| \frac{1}{\mu_t^n(\mathbb{R}^{d_1})} \int_{\mathbb{R}^{d_1}} \varphi(x) \mu_t^n(dx) - E^0\{\varphi(X_t) | \{Y_{t_k}, t_k \le t\}\} \right| \to 0 \tag{2}$$

in probability, as $\varepsilon \to 0$, $n \to \infty$, with a rate of convergence for all continuous bounded $\varphi$ so long as $\inf_{\varepsilon,n}\{\varepsilon^{1/2} n\} > 0$. Indeed, we establish much more in terms of estimates on the error in (2) and types of convergence, including 2nd-mean and almost sure.

**2. Notation, results and algorithm.** In the current section we set our main notation, state our results and give our particle system algorithm to asymptotically solve our filtering problem. The proofs of the stated results are given in a later section. During the course of a proof we use the same symbol $c$ for constants, although the exact value of the constant may change. We show the dependence of $c$ on relevant parameters unless suppression causes no confusion. Throughout this note, we take $|\cdot|$ to be both Euclidean distance, as well as absolute value (depending on context). We fix a constant $T > 0$ and let $0 < \varepsilon \le 1$. To conserve space, we define

$$\langle \lambda, \varphi \rangle \doteq \int_{\mathbb{R}^{d_1}} \varphi(x) \lambda(dx)$$

for all signed Borel measures $\lambda$ and $|\lambda|$-integrable functions $\varphi$. Next, we let $\mathcal{B}_b(\mathbb{R}^{d_1})$ denote the set of all measurable bounded functions on $\mathbb{R}^{d_1}$. For $\varphi \in \mathcal{B}_b(\mathbb{R}^{d_1})$, we let $\|\varphi\|_\infty$ denote its supremum norm. We denote by $\mathcal{L}$ the generator of the signal $X$, and define

$$\mathcal{T} \doteq \{\varphi \in \mathcal{B}_b(\mathbb{R}^{d_1}) : \mathcal{L}\varphi \in \mathcal{B}_b(\mathbb{R}^{d_1})\}.$$

Then, one can check that $\mathcal{T}$ contains all finite multivariate trigonometric series. Further, we let $\mathcal{S}(\mathbb{R}^{d_1})$ denote the set of all rapidly decreasing functions on $\mathbb{R}^{d_1}$ and assume that $h = (h_1, \dots, h_{d_2})'$ with $h_i \in \mathcal{S}(\mathbb{R}^{d_1})$ for each $1 \le i \le d_2$. Finally, we let $\lfloor u \rfloor$ denote the greatest integer not more than a real number $u$, let $\lceil u \rceil$ denote the least integer not less than $u$, and adopt the convention that a product over zero or a negative number of elements is one.

We define filtration

$$\mathcal{Y}^t \doteq \sigma\{Y_{t_k}, t_k \le t\} \vee \mathcal{N}$$

for the observations $Y$, where $\mathcal{N}$ is the collection of $P^0$-null sets of $(\Omega, \mathcal{F})$. Motivated by the reference probability measure method for filtering, we define a new probability measure via

$$\frac{dP}{dP^0} \doteq \eta_T,$$

where

$$\eta_t \doteq \prod_{k=1}^{\lfloor t/\varepsilon \rfloor} \exp\biggl\{-h'(X_{t_k})(V_{t_k} - V_{t_{k-1}}) - \frac{(h'h)(X_{t_k})(t_k - t_{k-1})}{2}\biggr\} \tag{3}$$

for $0 \le t \le T$. We define $\mathcal{X}^t \doteq \sigma\{X_s, 0 \le s \le t\} \vee \mathcal{N}$ and find that $\{\eta_t, t \in [0,T]\}$ is an $\{\mathcal{X}^T \vee \mathcal{Y}^t\}_{0 \le t \le T}$-martingale with respect to $P^0$. Under $P$, $\{Y_{t_k} - Y_{t_{k-1}}, k = 1, 2, \ldots, \lfloor T/\varepsilon \rfloor\}$ is a sequence of $N(0, \varepsilon I_{d_2})$ random vectors independent of $X$ and the law of $X$ remains unchanged. Yet, by (3), it follows that

$$\eta_T^{-1} = \prod_{k=1}^{\lfloor T/\varepsilon \rfloor} \exp\biggl\{h'(X_{t_k})(Y_{t_k} - Y_{t_{k-1}}) - \frac{(h'h)(X_{t_k})(t_k - t_{k-1})}{2}\biggr\}.$$

We let $E$ be the expectation with respect to $P$ and define

$$\langle \mu_t, \varphi \rangle = E\{\varphi(X_t)\eta_t^{-1} | \mathcal{Y}^t\}$$

for $0 \le t \le T$. Then, it follows from Bayes' rule that, for any $\varphi \in \mathcal{B}_b(\mathbb{R}^{d_1})$,

$$\begin{aligned} E^0\{\varphi(X_t)|\mathcal{Y}^t\} &= \frac{E\{\varphi(X_t)\eta_T^{-1}|\mathcal{Y}^t\}}{E\{\eta_T^{-1}|\mathcal{Y}^t\}} \\ &= \frac{E\{\varphi(X_t)\eta_t^{-1}|\mathcal{Y}^t\}}{E\{\eta_t^{-1}|\mathcal{Y}^t\}} \\ &= \frac{\langle \mu_t, \varphi \rangle}{\langle \mu_t, 1 \rangle} \end{aligned}$$

by the $\mathcal{X}^T \vee \mathcal{Y}^t$-martingale property of $\eta_t^{-1}$ with respect to $P$. For the processes that we will work with later, one may always assume that $X$ is cádlág and, hence, that $\mu_t$ is also (cf. [19]). We always work with this cádlág version.

First, considering the optimal solution to the filtering problem, we have the following lemma whose proof is sketched in the Appendix.

LEMMA 1. *Suppose that $\mu_0$ is the distribution of the initial signal state. Then, $\{\mu_t, t \ge 0\}$ is the unique measure-valued, $\{\mathcal{Y}^t\}_{t \ge 0}$-progressive process satisfying*

$$\langle \mu_t, \varphi \rangle = \langle \mu_0, \varphi \rangle + \int_0^t \langle \mu_s, \mathcal{L}\varphi \rangle \, ds$$

$$+\sum_{k=1}^{\lfloor t/\varepsilon\rfloor}\left\langle \mu_{t_k-},\varphi\left(\exp\left\{(Y_{t_k}-Y_{t_{k-1}})'h-\frac{(t_k-t_{k-1})h'h}{2}\right\}-1\right)\right\rangle \tag{4}$$

*for all* $\varphi\in\mathcal{T}$.

Moving to our particle approximation, we recall that $\|\mu_t^n\|$ and $\mu_t^n$ denote, respectively, the number of particles alive and the empirical measure for the particles as in the Introduction. Once we have particles $\{X_t^{i,n}\}_{i=1}^{\|\mu_t^n\|}, t\geq 0$, we can form our approximation to $\mu_t$ via empirical measure (1). Therefore, our pressing need is to find a good generation method for the particles. We suggest using the algorithm below to produce particles whose empirical measure is shown in the sequel to converge nicely to $\{\mu_t, t\geq 0\}$.

To ease the notation in what follows, we define

$$\varrho_k^\varepsilon(x)\doteq\exp\left\{(Y_{t_k}-Y_{t_{k-1}})'h(x)-\frac{(t_k-t_{k-1})(h'h)(x)}{2}\right\}-1 \tag{5}$$

and

$$DY_t^\varepsilon(x)\doteq\sum_{k=1}^{\infty}\delta_{k\varepsilon}(t)\varrho_k^\varepsilon(x)$$

$$\xi_k^\varepsilon(x)\doteq\begin{cases}\varrho_k^\varepsilon(x), & \text{if } \varrho_k^\varepsilon(x)<0,\\ \varrho_k^\varepsilon(x)-\lfloor\varrho_k^\varepsilon(x)\rfloor, & \text{otherwise.}\end{cases}$$

Moreover, due to the fact that we have both continuous and discrete components to our systems, it will be convenient in the sequel to interpret $\delta_{k\varepsilon}$ in two ways:

$$\int_s^u\delta_{k\varepsilon}(t)\,dt=\begin{cases}1, & \text{if } k\varepsilon\in(s,u],\\ 0, & \text{otherwise,}\end{cases}$$

and

$$\sum_{j=i}^{l}\delta_{k\varepsilon}(j\varepsilon)=\begin{cases}1, & \text{if } k\in\{i,i+1,\dots,l\},\\ 0, & \text{otherwise.}\end{cases}$$

Next, we let $\{\rho^i\}_{i=1}^n$ be $n$ independent $\mathbb{R}^{d_1}$-valued random variables with the distribution $\mu_0$, let $\{\tilde{X}^i\}_{i=1}^\infty$ be a sequence of independent $\mathbb{R}^{d_1}$-valued Lévy-stable processes with index $\alpha$ and spectral measure $\Gamma$, let $\{U^{i,k}\}_{i,k=1}^\infty$ be a sequence of independent uniform random variables on $[0,1]$. We suppose that $\{\rho^i\}_{i=1}^n$, $\{\tilde{X}^i\}_{i=1}^\infty$ and $\{U^{i,k}\}_{i,k=1}^\infty$ are defined on the same probability space $(\Omega^*,\mathcal{F}^*,P^*)$ and they are independent of one another. We define the product probability space $(\hat{\Omega},\hat{\mathcal{F}},\hat{P})\doteq(\Omega\otimes\Omega^*,\mathcal{F}\otimes\mathcal{F}^*,P^0\otimes P^*)$ and let $\hat{E}$ be the expectation with respect to $\hat{P}$. Then, to construct our particle system to approximate $\mu_t$, we do the following:

1. Let $\|\mu_0^n\| = n$ and $X_0^{i,n} \doteq \rho^i$ for $i = 1, \dots, \|\mu_0^n\|$ (*Assign initial particle locations*).
2. For $k = 1, 2, \dots,$ do the following:
   (a) Set $X_t^{i,n} = X_{t_{k-1}}^{i,n} + (\tilde{X}_t^i - \tilde{X}_{t_{k-1}}^i)$ for $t \in [t_{k-1}, t_k)$ and $i \in \{1, \dots, \|\mu_{t_{k-1}}^n\|\}$ (*Evolve particles as signal*).
   (b) For $i = 1, \dots, \|\mu_{t_{k-1}}^n\|$, do
      (i) If $\varrho_k^\varepsilon(X_{t_k-}^{(i,n)}) \geq 0$ (*Branch*),
         (A) Replace particle $X_{t_k-}^{i,n}$ with $m \doteq \lfloor \varrho_k^\varepsilon(X_{t_k-}^{i,n}) \rfloor + 1$ particles $X_{t_k}^{(i,1),n}$, $\dots$, $X_{t_k}^{(i,m),n}$ at site $X_{t_k-}^{i,n}$,
         (B) Add 1 more particle $X_{t_k}^{(i,m+1),n}$ at site $X_{t_k-}^{i,n}$ if $U^{i,k} \leq \varrho_k^\varepsilon(X_{t_k-}^{i,n}) - \lfloor \varrho_k^\varepsilon(X_{t_k-}^{i,n}) \rfloor$,
      (ii) Otherwise,
         (A) Make no change if $U^{i,k} > |\varrho_k^\varepsilon(X_{t_k-}^{i,n})|$,
         (B) Kill $X_{t_k-}^{i,n}$ if $U^{i,k} \leq |\varrho_k^\varepsilon(X_{t_k-}^{i,n})|$ (*Particle will just be removed*).
3. Relabel the alive particles to be $\{X_{t_k}^{i,n}\}_{i=1}^{\|\mu_{t_k}^n\|}$ so that $\|\mu_{t_k}^n\|$ is the number of particles alive.

Our main contributions can be considered as the popularization of this algorithm and its analysis. As we already mentioned, $U^{i,k} \leq |\xi_k^\varepsilon(X_{t_k-}^{i,n})|$, hence, branching or killing will seldom occur at a particular observation for small $\varepsilon > 0$. In preparation to listing our main analytic results, we wish now to assert that our empirical measures or particle density profiles

$$\mu_t^n \doteq \frac{1}{n} \sum_{i=1}^{\|\mu_t^n\|} \delta_{X_t^{i,n}}$$

do henceforth pertain only to the particles $\{X_t^{i,n}\}_{i=1}^{\|\mu_t^n\|}$, $t \geq 0$, generated by this algorithm. We define new filtrations $\{\mathcal{F}^t\}_{t\geq 0}$, $\{\mathcal{G}^t\}_{t\geq 0}$ to keep track of current information in our empirical measures and our whole particle system construction via

$$\mathcal{F}^t \doteq \bigcap_{\delta>0} \sigma\{X_s^{i,n}, i = 1, \dots, \|\mu_s^n\|, s \leq t + \delta\} \vee \mathcal{Y}^t,$$

$$\mathcal{G}^t \doteq \bigcap_{\delta>0} \sigma\{X_s^{i,n}, i = 1, \dots, \|\mu_s^n\|, s \leq t + \delta\} \vee \mathcal{Y}^T \vee \sigma\{U^{i,k}, t_k \leq t, i = 1, 2, \dots\}.$$

Further, we interpret our particle system approximation as a (purely atomic) measure-valued cádlág process through the stochastic equation (6) in Proposition 2. Hereafter, for a semimartingale $Z$, we use $[Z] = \{[Z, Z]_t, t \geq 0\}$ to denote its quadratic variation process.

PROPOSITION 2. *Suppose that* $\{\mu_t^n, t \geq 0\}$ *is the particle density profile constructed by the preceding algorithm. Then*

$$\begin{aligned}\langle \mu_t^n, \varphi\rangle &= \langle \mu_0^n, \varphi\rangle + \int_0^t \langle \mu_s^n, \mathcal{L}\varphi\rangle\, ds \\ &\quad + \sum_{k=1}^{\lfloor t/\varepsilon\rfloor} \langle \mu_{k\varepsilon-}^n, \varrho_k^\varepsilon \varphi\rangle + \mathcal{M}_t^n(\varphi)\end{aligned} \tag{6}$$

*for all* $\varphi \in \mathcal{T}$, *where* $\{\mathcal{M}_t^n(\varphi)\}_{t\geq 0}$ *is a cádlág* $\{\mathcal{G}^t\}_{t\geq 0}$*-martingale. We define* $\hat{E}^U$ *to be the expectation taken only with respect to the* $\{U^{i,k}\}$. *Then*

$$\begin{aligned}\hat{E}^U\{[\mathcal{M}^n(\varphi)]_t\} &= \frac{1}{n^2}\sum_{k=0}^{\lfloor t/\varepsilon\rfloor}\sum_{i=1}^{\|\mu_{t_k}^n\|}([\varphi(X^{i,n})]_{t_{k+1}\wedge t} - [\varphi(X^{i,n})]_{t_k}) \\ &\quad + \frac{1}{n}\sum_{k=1}^{\lfloor t/\varepsilon\rfloor}\langle \mu_{k\varepsilon-}^n, (|\xi_k^\varepsilon| - (\xi_k^\varepsilon)^2)\varphi^2\rangle.\end{aligned} \tag{7}$$

*Moreover, we have that*

$$\begin{aligned}&\hat{E}\left\{\left|\sum_{k=\lfloor s/\varepsilon\rfloor+1}^{\lfloor t/\varepsilon\rfloor} f_k([\mathcal{M}^n(\varphi)]_{k\varepsilon} - \hat{E}^U\{[\mathcal{M}^n(\varphi)]_{k\varepsilon}\})\right|^r\right\} \\ &\qquad \leq c(r)(\|h'h\|_\infty \vee 1)^{r/2}\left|\sum_{k=\lfloor s/\varepsilon\rfloor+1}^{\lfloor t/\varepsilon\rfloor} f_k^2\right|^{r/2} \\ &\qquad\quad \times \frac{\varepsilon^{1/2}}{n^r}\left(\sup_{0\leq\tau\leq T}\hat{E}\{\langle\mu_\tau^n, 1\rangle^r\}\right)\|\varphi\|_\infty^{2r}\end{aligned} \tag{8}$$

*for any* $\{f_k\}_{k=1}^\infty \subset \mathbb{R}$, $0 \leq s < t \leq T$, *where* $r \geq 2$ *and* $c(r) > 0$ *is a constant independent of* $d_1$, $d_2$, $\varepsilon$, $n$, $t$, $s$, $\varphi$.

This representation lemma differs from standard formulations because it contains both continuous and discrete time components. It is possible to come up with a more complete martingale problem description by considering more general functionals $F(\langle \mu_t^n, \varphi\rangle)$ instead of just $\langle \mu_t^n, \varphi\rangle$. However, our representation is sufficient for our purposes. To prove Proposition 2, we need the following Lemma 3. The proofs of Lemma 3 and Proposition 2 are given in the Appendix.

LEMMA 3. *Suppose* $r \geq 1$. *Then, there is a constant* $c(r) > 0$ *independent of* $d_1$, $d_2$, $\varepsilon$, $x$, $k$ *such that*

$$\hat{E}\{|\varrho_k^\varepsilon(x)|^r\} \leq c(r)\|h'h\|_\infty^{r/2}\varepsilon^{r/2}$$

*for all* $x \in \mathbb{R}^{d_1}$ *and* $k = 1, 2, \ldots$.

By Lemma 1, $\{\mu_t, t \geq 0\}$ is the unique measure-valued, $\{\mathcal{Y}^t\}_{t\geq 0}$-progressive process such that

$$\langle \mu_t, \varphi\rangle = \langle \mu_0, \varphi\rangle + \int_0^t \langle \mu_{s-}, B_s^\varepsilon \varphi\rangle\, ds \tag{9}$$

for all $\varphi \in \mathcal{T}$, where

$$B_s^\varepsilon \varphi \doteq \mathcal{L}\varphi + DY_s^\varepsilon \varphi. \tag{10}$$

Note, here and in the sequel, integrals like $\int_0^t \langle \mu_{s-}, B_s^\varepsilon \varphi\rangle\, ds$ should be interpreted in the Lebesgue–Stieltjes sense, including jumps at $t$ but not at $0$ (owing to the fact that $DY_s^\varepsilon$ is a purely atomic measure and not a function). We let $\gamma < -d_1/2$ and define

$$\|\lambda\|_\gamma^2 \doteq \int_{\mathbb{R}^{d_1}} |\hat{\lambda}(\theta)|^2 \gamma(d\theta), \qquad \gamma(d\theta) \doteq (1+|\theta|^2)^\gamma\, d\theta,$$

$$\hat{\lambda}(\theta) \doteq \langle \lambda, e_{-\theta}\rangle, \qquad e_{-\theta}(x) \doteq e^{-i\theta' x} \qquad \forall\, \theta \in \mathbb{R}^{d_1},$$

where $\hat{\lambda}$ denotes the Fourier–Stieltjes transform for a signed measure $\lambda$. In the sequel, we use $\|\varphi\|_{L^2(\gamma)}$, $\|\varphi\|_{L^2}$, respectively, to denote the $L^2$-norm of a function $\varphi$ in $L^2(\mathbb{R}^{d_1}; \gamma(d\theta))$, $L^2(\mathbb{R}^{d_1}; d\theta)$. We denote $\|\Gamma\| = \Gamma(S_{d_1})$. For $m \in \mathbb{N}$, we define

$$\langle\langle h\rangle\rangle_m \doteq \sup_{1\leq i\leq d_2, |\tau|\leq m} \left\{ \left\| \left( \prod_{1\leq j\leq d_1} (|x_j|+1)\right) D^\tau h_i \right\|_\infty \right\}, \tag{11}$$

where $\tau = (\tau_1, \ldots, \tau_{d_1})$ with $\tau_j \in \mathbb{Z}_+$ is a multi-index, $|\tau| = \sum_{j=1}^{d_1} \tau_j$ and $D^\tau = \partial^{|\tau|}/(\partial x_1^{\tau_1} \cdots \partial x_{d_1}^{\tau_{d_1}})$.

Now, we can state our main results.

THEOREM 4. *Let $\{\mu_t^n, t \geq 0\}$ be our particle density profile as described above. Suppose $\gamma < -(d_1/2 + 2\alpha)$.*

(i) *Let $\Xi > 0$ be a constant. Then, there is a constant $c(\Xi, d_1, \alpha, \|\Gamma\|, \langle\langle h\rangle\rangle_{[d_1-2\gamma]+2}, T) > 0$ independent of $\varepsilon$, $n$, $t$, $s$ such that*

$$\begin{aligned} &\hat{E}^{1/2}\{\|(\mu_t^n - \mu_t) - (\mu_s^n - \mu_s)\|_\gamma^2\} \\ &\quad \leq \frac{c(\Xi, d_1, \alpha, \|\Gamma\|, \langle\langle h\rangle\rangle_{[d_1-2\gamma]+2}, T)}{\varepsilon^{1/8} n^{1/2}} \\ &\qquad \times \left\{ (t-s)^{1/4} + (t-s) \right. \\ &\qquad\quad \left. + \varepsilon^{1/2}\left(\left\lfloor \frac{t}{\varepsilon}\right\rfloor - \left\lfloor \frac{s}{\varepsilon}\right\rfloor\right)^{1/2} + \varepsilon^{1/4}\left(\left\lfloor \frac{t}{\varepsilon}\right\rfloor - \left\lfloor \frac{s}{\varepsilon}\right\rfloor\right)^{1/4} \right\} \\ &\qquad\qquad\qquad \forall\, 0 \leq s < t \leq T \end{aligned} \tag{12}$$

*for any* $0 < \varepsilon \le 1$ *and* $n \in \mathbb{N}$ *satisfying* $\varepsilon^{1/2} n \ge \Xi$.

(ii) (Rate of convergence.) *Let* $0 < \varepsilon \le 1$ *be a constant. Then, there is a constant* $c(\varepsilon, d_1, \alpha, \|\Gamma\|, \langle\!\langle h \rangle\!\rangle_{[d_1 - 2\gamma]+2}, T) > 0$ *independent of* $n$ *such that*

$$\sup_{0 \le t \le T} \hat{E}^{1/2}\{\|\mu_t^n - \mu_t\|_\gamma^2\} \le \frac{c(\varepsilon, d_1, \alpha, \|\Gamma\|, \langle\!\langle h \rangle\!\rangle_{[d_1 - 2\gamma]+2}, T)}{n^{1/2}} \tag{13}$$

*for all* $n \in \mathbb{N}$.

COROLLARY 5. *Let* $\{\mu_t^n, t \ge 0\}$ *be our particle density profile as described above. Suppose that* $\alpha = 2$ *and* $\gamma < -(d_1/2 + 4)$.

(i) *Let* $\Xi > 0$ *and* $\beta > 1/8$ *be two constants. Then, there is a constant* $c(\Xi, d_1, \beta, \|\Gamma\|, \langle\!\langle h \rangle\!\rangle_{[d_1 - 2\gamma]+2}, T) > 0$ *independent of* $\varepsilon$, $n$ *such that*

$$\begin{aligned} &\hat{E}^{1/2}\left\{\sup_{0 \le s < t \le T} \|(\mu_t^n - \mu_t) - (\mu_s^n - \mu_s)\|_\gamma^2\right\} \\ &\qquad \le \frac{c(\Xi, d_1, \beta, \|\Gamma\|, \langle\!\langle h \rangle\!\rangle_{[d_1 - 2\gamma]+2}, T)}{\varepsilon^\beta n^{1/2}} \end{aligned} \tag{14}$$

*for any* $0 < \varepsilon \le 1$ *and* $n \in \mathbb{N}$ *satisfying* $\varepsilon^{1/2} n \ge \Xi$.

(ii) (Rate of convergence.) *Let* $0 < \varepsilon \le 1$ *be a constant. Then, there is a constant* $c(\varepsilon, d_1, \|\Gamma\|, \langle\!\langle h \rangle\!\rangle_{[d_1 - 2\gamma]+2}, T) > 0$ *independent of* $n$ *such that*

$$\hat{E}^{1/2}\left\{\sup_{0 \le t \le T} \|\mu_t^n - \mu_t\|_\gamma^2\right\} \le \frac{c(\varepsilon, d_1, \|\Gamma\|, \langle\!\langle h \rangle\!\rangle_{[d_1 - 2\gamma]+2}, T)}{n^{1/2}} \tag{15}$$

*for all* $n \in \mathbb{N}$.

REMARK 6. For the interacting mechanism chosen in the work of Del Moral [7], the number of particles remains constant and particles redistribute themselves around existing particle sites according to a multinomial distribution at observation times. Specifically, suppose $\{X_{t_k-}^{1,n}, \dots, X_{t_k-}^{n,n}\}$ denotes the $n$ particle locations used to approximate the filtering problem solution just prior to $t_k$, $\{W_k^{1,n}, \dots, W_k^{n,n}\}$ are the normalized weights for the particles, and $\{X_{t_k}^{1,n}, \dots, X_{t_k}^{n,n}\}$ is the system immediately following the interaction. Then, the $X_{t_k}^{i,n}$'s are obtained from the $X_{t_k-}^{j,n}$'s by having each $X_{t_k}^{i,n}$ choose starting location $X_{t_k-}^{j,n}$ with probability $W_k^{j,n}$ independent of all other particle decisions. Each weight $W_k^{j,n}$ is a function of all the previous generation particles $\{X_{t_k-}^{1,n}, \dots, X_{t_k-}^{n,n}\}$, the current observation $Y_{t_k}$, the conditional distribution of the observation given the current signal state $X_{t_k}$, and the conditional distribution of signal $X_{t_k}$ given all the previous observations $\{Y_{t_j}, j < k\}$. Clearly, $\sum_{j=1}^n W_k^{j,n} = 1$ and the event $\{X_{t_k}^{i,n} \ne X_{t_k-}^{i,n}\}$ has probability $1 - W_k^{i,n}(\omega)$, so the expected number of branches or jumps created at an observation time is $n - 1$ even when the observation interval or the time between jumps is very small. Moreover, as mentioned in [5], the decision of where each particle will jump to requires sampling all particles, and the overall result is that a large amount of computational work must be done at observation times.

REMARK 7. In [4] rates of convergence for a branching particle approximation to the solution of the Zakai equation are deduced. For a class of test functions, exact rates of convergence are established for the filtering model with diffusion signal and continuous observations. The analysis in [4] hinges on a powerful representation formula of the variance of the branching mechanism in terms of the local time of an exponential martingale, which is quite different from the analysis in this paper. Throughout this paper Fourier analysis is used, which enables us to obtain powerful rates of convergence in Sobolev norms. (We refer the interested reader to [2] and references therein for some other works via Fourier analysis, which are close in spirit to our approach.) The analysis of the existing interacting and branching methods for continuous-discrete filters is rather complicated, as is evidenced by the limited number of existing estimates especially involving the time intervals between observations. As suggested in [4], the continuous observation time set-up makes the branching method converge slower. Our Theorem 4 and Corollary 5 reveal the subtle relationship between the number of initial particles and the length of the time intervals between observations. In particular, the convergence of the algorithm is ensured if $\inf_{\varepsilon,n}\{\varepsilon^{1/2}n\} > 0$. In a forthcoming work, we look forward to further developing the spectral method in this paper to obtain rates of convergence for more general (not necessary diffusion) Markov processes and other recently developed particle filters.

## 3. Proofs of Theorem 4 and Corollary 5.

### 3.1. *Auxiliary results used to establish Theorem 4 and Corollary 5.*

LEMMA 8. *Let $Z$ be a $\mathbb{R}^{d_1}$-valued Lévy-stable process on $(\hat{\Omega}, \hat{\mathcal{F}}, \hat{P})$ with index $\alpha \in (0,2]$ and spectral measure $\Gamma$. We define $\hat{Z}_t(\theta) \doteq e_{-\theta}(Z_t)$ and $\|[\hat{Z}(\theta)]_t\| \doteq [\operatorname{Re} \hat{Z}(\theta)]_t + [\operatorname{Im} \hat{Z}(\theta)]_t$, $\forall \theta \in \mathbb{R}^{d_1}$, $t \geq 0$. Then, for $0 \leq s < t < \infty$,*

$$\hat{E}\{\|[\hat{Z}(\theta)]_t\| - \|[\hat{Z}(\theta)]_s\|\} = 2(t-s)\int_{S_{d_1}} |\theta' z|^\alpha \Gamma(dz). \tag{16}$$

*Suppose $r > 1$. Then, there is a constant $c(r) > 0$ such that, for any $0 \leq s < t < \infty$,*

$$\hat{E}\{(\|[\hat{Z}(\theta)]_t\| - \|[\hat{Z}(\theta)]_s\|)^r\}$$

(17)
$$\leq (r)(t-s)\biggl\{\biggl(\int_{S_{d_1}} |\theta' z|^\alpha \Gamma(dz)\biggr) \vee \biggl(\int_{S_{d_1}} |\theta' z|^\alpha \Gamma(dz)\biggr)^r\biggr\}.$$

*Moreover, if* $\alpha = 2$*, then*

$$\|[\hat{Z}(\theta)]_t\| - \|[\hat{Z}(\theta)]_0\| = 2t \int_{S_{d_1}} |\theta' z|^2 \Gamma(dz). \tag{18}$$

PROOF. For $0 \leq s < t < \infty$, we let $\{\tau_j^m, j = 0, 1, \ldots, k_m\}_{m=1}^\infty$ be a refining sequence of partitions for $[s,t]$ with $s = \tau_0^m < \tau_1^m < \cdots < \tau_{k_m}^m = t$ and define

$$\pi_{s,t}^m \doteq \sum_{j=1}^{k_m} |\hat{Z}_{\tau_j^m}(\theta) - \hat{Z}_{\tau_{j-1}^m}(\theta)|^2,$$
$$\delta(\pi_{s,t}^m) \doteq \max_{1 \leq j \leq k_m} \{\tau_j^m - \tau_{j-1}^m\}.$$

Then, we find by direct calculation that

$$\begin{aligned}
&\hat{E}\{\|[\hat{Z}(\theta)]_t\| - \|[\hat{Z}(\theta)]_s\|\} \\
&\quad = \lim_{\delta(\pi_{s,t}^m) \to 0} \hat{E}\{\pi_{s,t}^m\} \\
&\quad = \lim_{\delta(\pi_{s,t}^m) \to 0} \hat{E}\Biggl\{\sum_{j=1}^{k_m} (2 - \{e_{-\theta}(Z_{\tau_j^m} - Z_{\tau_{j-1}^m}) + e_\theta(Z_{\tau_j^m} - Z_{\tau_{j-1}^m})\})\Biggr\} \\
&\quad = \lim_{\delta(\pi_{s,t}^m) \to 0} \begin{cases}
\displaystyle\sum_{j=1}^{k_m} 2\biggl(1 - \exp\biggl\{-(\tau_j^m - \tau_{j-1}^m) \int_{S_{d_1}} |\theta' z|^\alpha \Gamma(dz)\biggr\} \\
\qquad \times \cos\biggl((\tau_j^m - \tau_{j-1}^m) \displaystyle\int_{S_{d_1}} |\theta' z|^\alpha \operatorname{sign}(\theta' z) \\
\qquad\qquad \times \tan\biggl(\dfrac{\alpha\pi}{2}\biggr)\Gamma(dz)\biggr)\biggr), \\
\qquad\qquad\qquad \text{for } \alpha \neq 1, \\
\displaystyle\sum_{j=1}^{k_m} 2\biggl(1 - \exp\biggl\{-(\tau_j^m - \tau_{j-1}^m) \int_{S_{d_1}} |\theta' z| \Gamma(dz)\biggr\} \\
\qquad \times \cos\biggl((\tau_j^m - \tau_{j-1}^m) \displaystyle\int_{S_{d_1}} \frac{2}{\pi} |\theta' z| \operatorname{sign}(\theta' z) \\
\qquad\qquad \times \ln|\theta' z| \Gamma(dz)\biggr)\biggr), \\
\qquad\qquad\qquad \text{for } \alpha = 1
\end{cases} \\
&\quad = 2(t-s) \int_{S_{d_1}} |\theta' z|^\alpha \Gamma(dz).
\end{aligned}$$

By (16), to prove (17), we may assume without loss of generality that $r \in \mathbb{N}$. By the independence of the increments of $Z$, we find that

$$\begin{aligned}
&\hat{E}\{(\|[\hat{Z}(\theta)]_t\| - \|[\hat{Z}(\theta)]_s\|)^r\} \\
&\quad = \lim_{\delta(\pi^m_{s,t})\to 0} \hat{E}\{(\pi^m_{s,t})^r\} \\
&\quad = \lim_{\delta(\pi^m_{s,t})\to 0} \hat{E}\left\{\left(\sum_{j=1}^{k_m}(2 - \{e_{-\theta}(Z_{\tau_j^m} - Z_{\tau_{j-1}^m}) + e_\theta(Z_{\tau_j^m} - Z_{\tau_{j-1}^m})\})\right)^r\right\} \qquad (19)\\
&\quad = \lim_{\delta(\pi^m_{s,t})\to 0} \sum_{\substack{\alpha_1+\cdots+\alpha_{k_m}=r\\ \alpha_1,\ldots,\alpha_{k_m}\in\mathbb{Z}_+}} \binom{r}{\alpha_1,\ldots,\alpha_{k_m}} \\
&\qquad\quad \times \prod_{j=1}^{k_m} \hat{E}\{(2 - \{e_{-\theta}(Z_{\tau_j^m} - Z_{\tau_{j-1}^m}) + e_\theta(Z_{\tau_j^m} - Z_{\tau_{j-1}^m})\})^{\alpha_j}\}.
\end{aligned}$$

Note that, for $\alpha_j \geq 1,\ 1 \leq j \leq k_m$,

$$\begin{aligned}
&\hat{E}\{(2 - \{e_{-\theta}(Z_{\tau_j^m} - Z_{\tau_{j-1}^m}) + e_\theta(Z_{\tau_j^m} - Z_{\tau_{j-1}^m})\})^{\alpha_j}\} \\
&\quad = \sum_{l=0}^{\alpha_j}\left\{\binom{\alpha_j}{l} 2^l (-1)^{\alpha_j - l}\right. \\
&\qquad\quad \times \sum_{q=0}^{\alpha_j - l}\left(\binom{\alpha_j - l}{q}\right. \\
&\qquad\qquad \times \exp\left\{-(\tau_j^m - \tau_{j-1}^m)\right. \\
&\qquad\qquad \times \begin{cases}
\displaystyle\int_{S_{d_1}} |(2q + l - \alpha_j)\theta' z|^\alpha \\
\qquad \times \left(1 - i\,\mathrm{sign}((2q + l - \alpha_j)\theta' z)\right. \\
\qquad\qquad \left.\left.\left.\left.\times \tan\dfrac{\alpha\pi}{2}\right)\Gamma(dz)\right\}\right)\right\}, & \text{for } \alpha \neq 1, \\
\displaystyle\int_{S_{d_1}} |(2q + l - \alpha_j)\theta' z| \\
\qquad \times \left(1 + \dfrac{2i}{\pi}\mathrm{sign}((2q + l - \alpha_j)\theta' z)\right. \\
\qquad\qquad \left.\left.\left.\left.\times \ln|\theta' z|\right)\Gamma(dz)\right\}\right)\right\}, & \text{for } \alpha = 1
\end{cases} \qquad (20)
\end{aligned}$$

$$
\begin{aligned}
&= -\sum_{l=0}^{\alpha_j}\Bigg\{\binom{\alpha_j}{l}2^l(-1)^{\alpha_j-l} \\
&\qquad\times\sum_{q=0}^{\alpha_j-l}\Bigg(\binom{\alpha_j-l}{q}(\tau_j^m-\tau_{j-1}^m) \\
&\qquad\qquad\times\int_{S_{d_1}}|(2q+l-\alpha_j)\theta'z|^\alpha\Gamma(dz)\Bigg)\Bigg\}+O((\tau_j^m-\tau_{j-1}^m)^2) \\
&\le c(r)(\tau_j^m-\tau_{j-1}^m)\int_{S_{d_1}}|\theta'z|^\alpha\Gamma(dz)+O((\tau_j^m-\tau_{j-1}^m)^2).
\end{aligned}
$$

Thus, by(19) and (20), we find that

$$
\begin{aligned}
&\hat{E}\{(\|[\hat{Z}(\theta)]_t\|-\|[\hat{Z}(\theta)]_s\|)^r\} \\
&\quad\le \lim_{\delta(\pi_{s,t}^m)\to 0}\sum_{\substack{\alpha_1+\cdots+\alpha_{k_m}=r\\ \alpha_1,\ldots,\alpha_{k_m}\in\mathbb{Z}_+}}\binom{r}{\alpha_1,\ldots,\alpha_{k_m}} \\
&\qquad\times\prod_{\alpha_j\ge 1}\Bigg(c(r)(\tau_j^m-\tau_{j-1}^m)\int_{S_{d_1}}|\theta'z|^\alpha\Gamma(dz) \\
&\qquad\qquad+O((\tau_j^m-\tau_{j-1}^m)^2)\Bigg) \\
&\quad\le c(r)(t-s)\Bigg\{\Bigg(\int_{S_{d_1}}|\theta'z|^\alpha\Gamma(dz)\Bigg)\vee\Bigg(\int_{S_{d_1}}|\theta'z|^\alpha\Gamma(dz)\Bigg)^r\Bigg\}.
\end{aligned}
$$

If $\alpha=2$, then we find by the independence of the increments of $Z$ that

$$
\begin{aligned}
&\hat{E}\Bigg\{\Bigg(\|[\hat{Z}(\theta)]_t\|-\|[\hat{Z}(\theta)]_0\|-2t\int_{S_{d_1}}|\theta'z|^2\Gamma(dz)\Bigg)^2\Bigg\} \\
&\quad=\lim_{\delta(\pi_{0,t}^m)\to 0}\hat{E}\Bigg\{\Bigg(\pi_{0,t}^m-\sum_{j=1}^{k_m}2\Bigg(1-\exp\Bigg\{-(\tau_j^m-\tau_{j-1}^m) \\
&\qquad\qquad\times\int_{S_{d_1}}|\theta'z|^2\Gamma(dz)\Bigg\}\Bigg)\Bigg)^2\Bigg\} \\
&\quad=\lim_{\delta(\pi_{0,t}^m)\to 0}\sum_{j=1}^{k_m}\hat{E}\Bigg\{\Bigg(2\exp\Bigg\{-(\tau_j^m-\tau_{j-1}^m)\int_{S_{d_1}}|\theta'z|^2\Gamma(dz)\Bigg\} \\
&\qquad\qquad-\{e_{-\theta}(Z_{\tau_j^m}-Z_{\tau_{j-1}^m})+e_\theta(Z_{\tau_j^m}-Z_{\tau_{j-1}^m})\}\Bigg)^2\Bigg\}
\end{aligned}
$$

$$
\begin{aligned}
&= \lim_{\delta(\pi_{0,t}^m)\to 0} \sum_{j=1}^{k_m} \biggl(4\exp\biggl\{-2(\tau_j^m-\tau_{j-1}^m)\int_{S_{d_1}} |\theta' z|^2\Gamma(dz)\biggr\} \\
&\qquad\qquad - 8\exp\biggl\{-2(\tau_j^m-\tau_{j-1}^m)\int_{S_{d_1}} |\theta' z|^2\Gamma(dz)\biggr\} \\
&\qquad\qquad + 2\biggl(1+\exp\biggl\{-4(\tau_j^m-\tau_{j-1}^m)\int_{S_{d_1}} |\theta' z|^2\Gamma(dz)\biggr\}\biggr)\biggr) \\
&= 0.
\end{aligned}
$$

Therefore, (18) follows. □

Lemma 9. *Suppose that $r \ge 1$ and $\Xi > 0$ is a constant. Then, there is a constant $c(r,\Xi,T) > 0$ independent of $d_1$, $d_2$, $\varepsilon$, $n$ such that the empirical measure of our particle system satisfies*

$$\sup_{0\le t\le T} \hat{E}^{1/r}\{\langle \mu_t^n, 1\rangle^r\} \le c(r,\Xi,T)(\|h'h\|_\infty \vee 1)^r$$

*for any $0 < \varepsilon \le 1$ and $n \in \mathbb{N}$ satisfying $\varepsilon^{1/2} n \ge \Xi$.*

Proof. By (6), (59) in the Appendix, Lemma 3 and induction, one finds that $\sup_{0\le t\le T} \hat{E}\{\langle \mu_t^n, 1\rangle^r\} < \infty$. We define

$$\zeta_k^\varepsilon \doteq |\xi_k^\varepsilon| - (\xi_k^\varepsilon)^2. \tag{21}$$

From (6)–(8) with $\varphi = 1$, noting that $\{\sum_{k=1}^{\lfloor t/\varepsilon\rfloor} \langle \mu_{k\varepsilon-}^n, \varrho_k^\varepsilon\rangle\}_{t\ge 0}$ is an $\{\mathcal{F}_{t-}\}_{t\ge 0}$-martingale and using Burkholder's inequality, independence, Jensen's inequality, Lemma 3, Minkowski's integral inequality and (21), we find that

$$
\begin{aligned}
&\hat{E}\{\langle \mu_t^n, 1\rangle^r\} \\
&\quad \le c(r)\biggl\{\hat{E}\{\langle \mu_0^n, 1\rangle^r\} + \hat{E}\biggl\{\biggl(\sum_{k=1}^{\lfloor t/\varepsilon\rfloor} \langle \mu_{k\varepsilon-}^n, \varrho_k^\varepsilon\rangle\biggr)^r\biggr\} \\
&\qquad\qquad + (\hat{E}\{|\hat{E}^U\{[\mathcal{M}^n(1)]_t\}|^{r/2}\} + \hat{E}\{|[\mathcal{M}^n(1)]_t - \hat{E}^U\{[\mathcal{M}^n(1)]_t\}|^{r/2}\})\biggr\} \\
&\quad = c(r)\biggl\{1 + \hat{E}\biggl\{\biggl(\sum_{k=1}^{\lfloor t/\varepsilon\rfloor} \langle \mu_{k\varepsilon-}^n, \varrho_k^\varepsilon\rangle\biggr)^r\biggr\} \\
&\qquad\qquad + \biggl(\frac{1}{n^{r/2}}\hat{E}\biggl\{\biggl(\sum_{k=1}^{\lfloor t/\varepsilon\rfloor} \langle \mu_{k\varepsilon-}^n, \zeta_k^\varepsilon\rangle\biggr)^{r/2}\biggr\} \\
&\qquad\qquad\quad + \hat{E}\{|[\mathcal{M}^n(1)]_t - \hat{E}^U\{[\mathcal{M}^n(1)]_t\}|^{r/2}\}\biggr)\biggr\}
\end{aligned}
$$

$$\le c(r)\Bigg\{1+\left\lfloor\frac{t}{\varepsilon}\right\rfloor^{(r/2)-1}\|h'h\|_\infty^{r/2}\varepsilon^{r/2}\sum_{k=1}^{\lfloor t/\varepsilon\rfloor}\hat{E}\{\langle\mu^n_{(k-1)\varepsilon},1\rangle^r\}$$

$$+\left\lfloor\frac{t}{\varepsilon}\right\rfloor^{(r/2)-1}\frac{\varepsilon^{r/4}}{n^{r/2}}\sup_{0\le s\le T}\hat{E}\{\langle\mu^n_s,1\rangle^{r/2}\}$$

$$+(\|h'h\|_\infty\vee 1)^{r/2}\left\lfloor\frac{t}{\varepsilon}\right\rfloor^{r/4}\frac{\varepsilon^{1/2}}{n^{r/2}}\sup_{0\le s\le T}\hat{E}\{\langle\mu^n_s,1\rangle^{r/2}\}\Bigg\},$$

where we have assumed without loss of generality that $r\ge 4$ above. Applying the discrete version of Gronwall's inequality, one thus discovers that

$$\sup_{0\le t\le T}\hat{E}\{\langle\mu^n_t,1\rangle^r\}$$

$$\le c(r,T)(\|h'h\|_\infty\vee 1)^{r/2}\left(1+\frac{\varepsilon^{1/2}}{\varepsilon^{r/4}n^{r/2}}\sup_{0\le t\le T}\hat{E}\{\langle\mu^n_t,1\rangle^{r/2}\}\right)$$

$$\le c(r,T)(\|h'h\|_\infty\vee 1)^{r/2}\left(1+\frac{\varepsilon^{1/2}}{\varepsilon^{r/4}n^{r/2}}\left(\sup_{0\le t\le T}\hat{E}\{\langle\mu^n_t,1\rangle^{r}\}\right)^{1/2}\right).$$

Therefore,

$$\sup_{0\le t\le T}\hat{E}\{\langle\mu^n_t,1\rangle^r\}$$

$$\le\Bigg(\Bigg(\frac{c(r,T)(\|h'h\|_\infty\vee 1)^{r/2}\varepsilon^{1/2}}{\varepsilon^{r/4}n^{r/2}}$$

$$+\sqrt{\frac{c^2(r,T)(\|h'h\|_\infty\vee 1)^{r}\varepsilon}{\varepsilon^{r/2}n^{r}}+4c(r,T)(\|h'h\|_\infty\vee 1)^{r/2}}\Bigg)\Big/2\Bigg)^2$$

$$\le c(r,\Xi,T)(\|h'h\|_\infty\vee 1)^{r}$$

for any $0<\varepsilon\le 1$ and $n\in\mathbb{N}$ satisfying $\varepsilon^{1/2}n\ge\Xi$. □

The following maximal inequality is a consequence of a theorem of Longnecker and Serfling [13] (cf. also [12]) and is used in (14) above.

Lemma 10. *Let $0\le U_1<U_2<\infty$ and suppose that $\{Q_t,U_1\le t\le U_2\}$ is a process assuming values in some normed vector space $(\mathcal{Z},\|\cdot\|)$ with the following conditions:* (i) *$t\to Q_t(\omega)$ is right continuous on $[U_1,U_2]$ for almost all $\omega$,* (ii) *There exist constants $\mu>1$ and $\nu>0$ such that $E\{\|Q_t-Q_s\|^\nu\}\le(h(s,t))^\mu$ for all $U_1\le s<t\le U_2$, where $h(t,s)$ is a nonnegative function satisfying $h(s,t)+h(t,u)\le h(s,u)$ for all $U_1\le s<t<u\le U_2$. Then, there exists a constant $A_{\mu,\nu}$ depending only upon $\mu$, $\nu$ such that*

$$E\Big\{\sup_{U_1\le s<t\le U_2}\|Q_t-Q_s\|^\nu\Big\}\le A_{\mu,\nu}(h(U_1,U_2))^\mu.$$

PROOF. Let $\{t_i^l, i=0,1,\dots,n_l\}_{l=1}^{\infty}$ be a refining sequence of partitions for $[U_1,U_2]$ with $U_1 = t_0^l < t_1^l < \cdots < t_{n_l}^l = U_2$ and define

$$\tau_k^l \doteq Q_{t_k^l} - Q_{t_{k-1}^l}, \qquad g_l(i,j) \doteq h(t_j^l, t_{i-1}^l) \qquad \forall\, i,j,k \in \{1,\dots,n_l\}, i<j.$$

Then, we can apply Theorem 1 of [13] to find that there is a constant $A_{\mu,\nu}$ depending only upon $\mu$, $\nu$ such that

$$\begin{aligned} E\Big\{\sup_{U_1\le t_i^l<t_j^l\le U_2} \|Q_{t_j^l} - Q_{t_i^l}\|^\nu\Big\} &= E\Bigg\{\sup_{1\le i<j\le n_l}\Bigg\|\sum_{k=i+1}^{j}\tau_k^l\Bigg\|^\nu\Bigg\} \\ &\le 2^\nu\Bigg(\sup_{1\le j\le n_l}\Bigg\|\sum_{k=1}^{j}\tau_k\Bigg\|^\nu\Bigg) \\ &\le A_{\mu,\nu}(g_l(1,n_l))^\mu \\ &= A_{\mu,\nu}(h(U_1,U_2))^\mu. \end{aligned}$$

The lemma therefore follows from monotone convergence and the observation that right continuity guarantees that

$$\sup_{U_1\le t_i^l<t_j^l\le U_2}\|Q_{t_j^l}-Q_{t_i^l}\| \overset{l\to\infty}{\nearrow} \sup_{U_1\le s<t\le U_2}\|Q_t - Q_s\|^\nu. \qquad \square$$

3.2. *Proof of Theorem 4.* Recalling (5), (6), (9) and (10), we find that $\mu^n - \mu$ satisfies

$$\langle \mu_t^n - \mu_t, \varphi\rangle = \langle \mu_0^n - \mu_0, \varphi\rangle + \int_0^t \langle \mu_{s-}^n - \mu_{s-}, B_s^\varepsilon \varphi\rangle\, ds + \mathcal{M}_t^n(\varphi)$$

for all $\varphi \in \mathcal{T}$, where $\mathcal{M}_t^n(\varphi)$ is the martingale of Proposition 2. We define

$$\ell(\theta) \doteq \begin{cases} -\displaystyle\int_{S_{d_1}} |\theta' z|^\alpha \left(1 + i\,\mathrm{sign}(\theta' z)\tan\left(\frac{\alpha\pi}{2}\right)\right)\Gamma(dz), & \text{for } \alpha \neq 1, \\ -\displaystyle\int_{S_{d_1}} |\theta' z|\left(1 - \frac{2i}{\pi}\,\mathrm{sign}(\theta' z)\ln|\theta' z|\right)\Gamma(dz), & \text{for } \alpha = 1. \end{cases}$$

Then, using $\varphi = e_{-\theta}$, we find that

$$\begin{aligned} &\langle \mu_t^n - \mu_t, e_{-\theta}\rangle \\ &\quad = \langle \mu_0^n - \mu_0, e_{-\theta}\rangle \\ &\qquad + \int_0^t \langle \mu_{s-}^n - \mu_{s-}, \ell(\theta)e_{-\theta} + DY_s^\varepsilon e_{-\theta}\rangle\, ds + \hat{\mathcal{M}}_t^n(\theta) \\ &\qquad\qquad \forall\, \theta \in \mathbb{R}^{d_1}. \end{aligned} \tag{22}$$

Hereafter, to ease the notation, we let $\hat{\mathcal{M}}_t^n(\theta) = \mathcal{M}_t^n(e_{-\theta})$. We define

$$\|[\hat{\mathcal{M}}^n(\theta)]_t\| \doteq [\mathrm{Re}\,\hat{\mathcal{M}}^n(\theta)]_t + [\mathrm{Im}\,\hat{\mathcal{M}}^n(\theta)]_t,$$
$$\|[\hat{X}^{i,n}(\theta)]_t\| \doteq [\mathrm{Re}\,\hat{X}^{i,n}(\theta)]_t + [\mathrm{Im}\,\hat{X}^{i,n}(\theta)]_t, \qquad \hat{X}_t^{i,n}(\theta) \doteq e_{-\theta}(X_t^{i,n}).$$

Then, from Proposition 2 and (21), we find that $\{\hat{\mathcal{M}}_t^n(\theta)\}_{t\ge 0}$ is a complex martingale with

$$\hat{E}^U\{\|[\hat{\mathcal{M}}^n(\theta)]_t\|\} = \frac{1}{n^2}\sum_{k=0}^{\lfloor t/\varepsilon\rfloor}\sum_{i=1}^{\|\mu_{t_k}^n\|}(\|[\hat{X}^{i,n}(\theta)]_{t_{k+1}\wedge t}\| - \|[\hat{X}^{i,n}(\theta)]_{t_k}\|) + \frac{1}{n}\sum_{k=1}^{\lfloor t/\varepsilon\rfloor}\langle \mu_{k\varepsilon-}^n, \zeta_k^\varepsilon\rangle. \tag{23}$$

Next, we divide $\langle \mu_t^n - \mu_t, e_{-\theta}\rangle$ into components:

$$\langle \mu_t^n - \mu_t, e_{-\theta}\rangle = \hat{u}_t^n(\theta) + \hat{v}_t^n(\theta) + \hat{\chi}_t^n(\theta).$$

Here, we define

$$\hat{u}_t^n(\theta) \doteq \int_0^t \ell(\theta)\hat{u}_s^n(\theta)\,ds + \hat{\mathcal{M}}_t^n(\theta), \tag{24}$$

$$\hat{\chi}_t^n(\theta) \doteq \langle \chi_t^n, e_{-\theta}\rangle \tag{25}$$

with

$$\langle \chi_t^n, \varphi\rangle \doteq \langle \mu_0^n - \mu_0, \varphi\rangle + \int_0^t \langle \chi_{s_-}^n, B_s^\varepsilon\varphi\rangle\,ds \qquad \forall\,\varphi\in\mathcal{T}, \tag{26}$$

and

$$\hat{v}_t^n(\theta) \doteq \langle \mu_t^n - \mu_t, e_{-\theta}\rangle - \hat{u}_t^n(\theta) - \hat{\chi}_t^n(\theta). \tag{27}$$

Note that, in the above definition, $\chi_t^n$ is just the unnormalized filtering process $\mu_t$ with the initial distribution $\mu_0^n - \mu_0$. We define

$$A_1 \doteq \hat{E}^{1/2}\{\|\hat{u}_t^n - \hat{u}_s^n\|_{L^2(\gamma)}^2\},$$
$$A_2 \doteq \hat{E}^{1/2}\{\|\hat{v}_t^n - \hat{v}_s^n\|_{L^2(\gamma)}^2\},$$
$$A_3 \doteq \hat{E}^{1/2}\{\|\hat{\chi}_t^n - \hat{\chi}_s^n\|_{L^2(\gamma)}^2\}.$$

Then,

$$E^{1/2}\{\|(\mu_t^n - \mu_t) - (\mu_s^n - \mu_s)\|_\gamma^2\} \le A_1 + A_2 + A_3 \tag{28}$$

by Minkowski's inequality. In the following, we will estimate $A_i$, $1\le i\le 3$, one by one.

(a) *Estimation of $A_1$.*

One finds from Proposition 2 that the following Wiener integral makes sense and from (24), as well as integration by parts, that

$$\hat{u}_t^n(\theta) = \int_0^t \exp\{(t-s)\ell(\theta)\}\, d\hat{\mathcal{M}}_s^n(\theta). \tag{29}$$

Fixing a $r \geq 2$, one finds from (29) that, for $0 \leq s < t \leq T$,

$$\begin{aligned}
&\hat{E}\{|\hat{u}_t^n(\theta) - \hat{u}_s^n(\theta)|^{2r}\} \\
&\quad = \hat{E}\left\{\left|(\exp\{(t-s)\ell(\theta)\} - 1)\hat{u}_s^n(\theta) + \int_s^t \exp\{(t-\tau)\ell(\theta)\}\, d\hat{\mathcal{M}}_\tau^n(\theta)\right|^{2r}\right\}.
\end{aligned}$$

Yet, using Burkholder's inequality, we find that

$$\begin{aligned}
&\hat{E}\{|\hat{u}_t^n(\theta) - \hat{u}_s^n(\theta)|^{2r}\} \\
&\quad \leq c(r)\Bigg\{|\exp\{(t-s)\ell(\theta)\} - 1|^{2r} \\
&\qquad\qquad \times \hat{E}\Bigg\{\Bigg(\int_0^s \exp\Bigg\{-2(s-\tau) \\
&\qquad\qquad\qquad \times \int_{S_{d_1}} |\theta' z|^\alpha \Gamma(dz)\Bigg\}\, d\|[\hat{\mathcal{M}}^n(\theta)]_\tau\|\Bigg)^r\Bigg\} \\
&\qquad\qquad + \hat{E}\Bigg\{\Bigg(\int_s^t \exp\Bigg\{-2(t-\tau) \\
&\qquad\qquad\qquad \times \int_{S_{d_1}} |\theta' z|^\alpha \Gamma(dz)\Bigg\}\, d\|[\hat{\mathcal{M}}^n(\theta)]_\tau\|\Bigg)^r\Bigg\}\Bigg\}.
\end{aligned} \tag{30}$$

We define

$$\mathcal{M}_\tau^{n,e}(\theta) \doteq \frac{1}{n}\sum_{k=0}^{\lfloor \tau/\varepsilon\rfloor}\sum_{i=1}^{\|\mu_{t_k}^n\|}\left(e_{-\theta}(X_{t_{k+1}\wedge\tau}^{i,n}) - e_{-\theta}(X_{t_k}^{i,n}) - \int_{t_k}^{t_{k+1}\wedge\tau}(\mathcal{L}e_{-\theta})(X_u^{i,n})\,du\right)$$

and

$$\begin{aligned}
\mathcal{M}_\tau^{n,b}(\theta) \doteq \frac{1}{n}\sum_{k=1}^{\lfloor \tau/\varepsilon\rfloor}\sum_{i=1}^{\|\mu_{t_k-}^n\|} &\langle \delta_{X_{t_k-}^{i,n}}, e_{-\theta}\rangle \\
&\times (\operatorname{sign}(\xi_k^\varepsilon(X_{t_k-}^{i,n}))1_{\{U^{i,k}\in[0,|\xi_k^\varepsilon(X_{t_k-}^{i,n})|)\}} - \xi_k^\varepsilon(X_{t_k-}^{i,n})).
\end{aligned}$$

Then, $\hat{\mathcal{M}}_\tau^{n,e}(\theta)$ and $\hat{\mathcal{M}}_\tau^{n,b}(\theta)$ are, respectively, the evolving and branching portions of the martingale $\hat{\mathcal{M}}_\tau^n(\theta)$. Considering (23) and separating $\hat{u}_t^n(\theta)$ into parts driven by $\hat{\mathcal{M}}_\tau^{n,e}(\theta)$ and $\hat{\mathcal{M}}_\tau^{n,b}(\theta)$, we find from double use of

Hölder's inequality and Lemma 8 that the evolving part of (30) satisfies

$$
\begin{aligned}
&\hat{E}^{1/r}\{|\hat{u}_t^{n,e}(\theta)-\hat{u}_s^{n,e}(\theta)|^{2r}\}\\
&\quad\le c(r)\Big\{|\exp\{(t-s)\ell(\theta)\}-1|^{2r}\hat{E}\{\|[\hat{\mathcal{M}}^{n,e}(\theta)]_s\|^r\}\\
&\qquad+\Big(\hat{E}\Big\{\int_s^t\exp\Big\{-4r(t-\tau)\int_{S_{d_1}}|\theta' z|^\alpha\Gamma(dz)\Big\}d\|[\hat{\mathcal{M}}^{n,e}(\theta)]_\tau\|\Big\}\Big)^{1/2} \qquad (31)\\
&\qquad\qquad\times(\hat{E}\{(\|[\hat{\mathcal{M}}^{n,e}(\theta)]_t\|-\|[\hat{\mathcal{M}}^{n,e}(\theta)]_s\|)^{2r-1}\})^{1/2}\Big\}^{1/r}\\
&\quad\le\frac{c(r,\|\Gamma\|)}{n}\sup_{0\le\tau\le T}\hat{E}^{1/r}\{\langle\mu_\tau^n,1\rangle^r\}\\
&\qquad\times\Big\{|\exp\{(t-s)\ell(\theta)\}-1|^{2r}(|\theta|^\alpha\vee|\theta|^{\alpha r})s\\
&\qquad\quad+\Big(1-\exp\Big\{-4r(t-s)\int_{S_{d_1}}|\theta' z|^\alpha\Gamma(dz)\Big\}\Big)^{1/2}\\
&\qquad\qquad\times((|\theta|^\alpha\vee|\theta|^{\alpha r})(t-s))^{1/2}\Big\}^{1/r}\\
&\quad\le\frac{c(r,\alpha,\|\Gamma\|)}{n}\sup_{0\le\tau\le T}\hat{E}^{1/r}\{\langle\mu_\tau^n,1\rangle^r\}\\
&\qquad\times\{(|\theta|^\alpha\vee|\theta|^{\alpha r})(|\theta|^\alpha|\ln|\theta\|+1)(t-s)\}^{1/r}.
\end{aligned}
$$

Furthermore, using the last two claims of Proposition 2 and (21), we find that the branching part of $\hat{u}_t^n(\theta)$ satisfies

$$
\begin{aligned}
&\hat{E}^{1/r}\{|\hat{u}_t^{n,b}(\theta)-\hat{u}_s^{n,b}(\theta)|^{2r}\}\\
&\quad\le\frac{c(r)}{n}\Big\{|\exp\{(t-s)\ell(\theta)\}-1|^{2r}\\
&\qquad\times\Big(\hat{E}\Big\{\Big|\sum_{k=1}^{\lfloor s/\varepsilon\rfloor}\exp\Big\{2(k\varepsilon-s)\int_{S_{d_1}}|\theta' z|^\alpha\Gamma(dz)\Big\}\langle\mu_{k\varepsilon-}^n,\zeta_k^\varepsilon\rangle\Big|^r\Big\}\\
&\qquad\quad+n^r\hat{E}\Big\{\Big|\sum_{k=1}^{\lfloor s/\varepsilon\rfloor}\exp\Big\{2(k\varepsilon-s)\int_{S_{d_1}}|\theta' z|^\alpha\Gamma(dz)\Big\}\\
&\qquad\qquad\times(\|[\hat{\mathcal{M}}^n(\theta)]_{k\varepsilon}\|-\hat{E}^U\{\|[\hat{\mathcal{M}}^n(\theta)]_{k\varepsilon}\|\})\Big|^r\Big\}\Big)
\end{aligned}
$$

$$
\begin{aligned}
&+\exp\left\{-2rt\int_{S_{d_1}}|\theta' z|^\alpha\Gamma(dz)\right\}\\
&\quad\times\left(\hat{E}\left\{\left|\sum_{k=\lfloor s/\varepsilon\rfloor+1}^{\lfloor t/\varepsilon\rfloor}\exp\left\{2k\varepsilon\int_{S_{d_1}}|\theta' z|^\alpha\Gamma(dz)\right\}\langle\mu^n_{k\varepsilon-},\zeta^\varepsilon_k\rangle\right|^r\right\}\right.\\
&\qquad+n^r\hat{E}\left\{\left|\sum_{k=\lfloor s/\varepsilon\rfloor+1}^{\lfloor t/\varepsilon\rfloor}\exp\left\{2k\varepsilon\int_{S_{d_1}}|\theta' z|^\alpha\Gamma(dz)\right\}\right.\right.\\
&\qquad\qquad\left.\left.\left.\times(\|[\hat{\mathcal{M}}^n(\theta)]_{k\varepsilon}\|-\hat{E}^U\{\|[\hat{\mathcal{M}}^n(\theta)]_{k\varepsilon}\|\})\right|^r\right\}\right)\Bigg\}^{1/r} \qquad (32)\\
&\le\frac{c(r)(\|h'h\|_\infty\vee 1)^{r/2}}{n}\\
&\quad\times\Bigg\{|\exp\{(t-s)\ell(\theta)\}-1|^{2r}\\
&\qquad\times\left(\hat{E}\left\{\left|\sum_{k=1}^{\lfloor s/\varepsilon\rfloor}\langle\mu^n_{k\varepsilon-},\zeta^\varepsilon_k\rangle\right|^r\right\}+\varepsilon^{(1-r)/2}\sup_{0\le\tau\le T}\hat{E}\{\langle\mu^n_\tau,1\rangle^r\}\right)\\
&\qquad+\hat{E}\left\{\left|\sum_{k=\lfloor s/\varepsilon\rfloor+1}^{\lfloor t\varepsilon\rfloor}\langle\mu^n_{k\varepsilon-},\zeta^\varepsilon_k\rangle\right|^r\right\}\\
&\qquad\qquad+\varepsilon^{1/2}\left(\left\lfloor\frac{t}{\varepsilon}\right\rfloor-\left\lfloor\frac{s}{\varepsilon}\right\rfloor\right)^{r/2}\sup_{0\le\tau\le T}\hat{E}\{\langle\mu^n_\tau,1\rangle^r\}\Bigg\}^{1/r}.
\end{aligned}
$$

Using Jensen's inequality applied to normalized sums and Lemma 3, we find from (32) that

$$
\begin{aligned}
&\hat{E}^{1/r}\{|\hat{u}^{n,b}_t(\theta)-\hat{u}^{n,b}_s(\theta)|^{2r}\}\\
&\quad\le\frac{c(r)(\|h'h\|_\infty\vee 1)^{r/2}\sup_{0\le\tau\le T}\hat{E}^{1/r}\{\langle\mu^n_\tau,1\rangle^r\}}{n}\\
&\qquad\times\left\{|\exp\{(t-s)\ell(\theta)\}-1|^{2r}\left(\varepsilon^{r/2}\left(\left\lfloor\frac{s}{\varepsilon}\right\rfloor\right)^{r-1}+\varepsilon^{(1-r)/2}\right)\right.\\
&\qquad\quad\left.+\varepsilon^{r/2}\left(\left\lfloor\frac{t}{\varepsilon}\right\rfloor-\left\lfloor\frac{s}{\varepsilon}\right\rfloor\right)^{r-1}+\varepsilon^{1/2}\left(\left\lfloor\frac{t}{\varepsilon}\right\rfloor-\left\lfloor\frac{s}{\varepsilon}\right\rfloor\right)^{r/2}\right\}^{1/r} \qquad (33)\\
&\quad\le\frac{c(r,\alpha,\|\Gamma\|)(\|h'h\|_\infty\vee 1)^{r/2}\sup_{0\le\tau\le T}\hat{E}^{1/r}\{\langle\mu^n_\tau,1\rangle^r\}}{\varepsilon^{1/2}n}
\end{aligned}
$$

$$\times\left\{\varepsilon^{1/(2r)}(1\wedge\{|\theta|^{\alpha}|\ln|\theta||(t-s)\})^2 + \varepsilon\left(\left\lfloor\frac{t}{\varepsilon}\right\rfloor-\left\lfloor\frac{s}{\varepsilon}\right\rfloor\right)^{(r-1)/r} + \varepsilon^{(r+1)/(2r)}\left(\left\lfloor\frac{t}{\varepsilon}\right\rfloor-\left\lfloor\frac{s}{\varepsilon}\right\rfloor\right)^{1/2}\right\}.$$

Piecing together (31), (33) and Lemma 9, one has that

$$\begin{aligned}\hat{E}^{1/r}&\{|\hat{u}_t^n(\theta)-\hat{u}_s^n(\theta)|^{2r}\}\\ &\le \frac{c(r,\Xi,\alpha,\|\Gamma\|,T)(\|h'h\|_\infty\vee 1)^{3r/2}}{\varepsilon^{1/2}n}\\ &\quad\times\left\{\varepsilon^{1/(2r)}(|\theta|^{\alpha/r}\vee|\theta|^{\alpha(r+2)/r})(t-s)^{1/r} + \varepsilon\left(\left\lfloor\frac{t}{\varepsilon}\right\rfloor-\left\lfloor\frac{s}{\varepsilon}\right\rfloor\right)^{(r-1)/r} + \varepsilon^{(r+1)/(2r)}\left(\left\lfloor\frac{t}{\varepsilon}\right\rfloor-\left\lfloor\frac{s}{\varepsilon}\right\rfloor\right)^{1/2}\right\}.\end{aligned}\tag{34}$$

Then, using Minkowski's integral inequality and (34), we find that

$$\begin{aligned}\hat{E}^{1/(2r)}&\{\|\hat{u}_t^n-\hat{u}_s^n\|_{L^2(\gamma)}^{2r}\}\\ &\le\left(\int_{\mathbb{R}^{d_1}}\hat{E}^{1/r}\{|\hat{u}_t^n(\theta)-\hat{u}_s^n(\theta)|^{2r}\}\gamma(d\theta)\right)^{1/2}\\ &\le\frac{c(r,\Xi,d_1,\alpha,\|\Gamma\|,\|h'h\|_\infty,T)}{\varepsilon^{1/4}n^{1/2}}\\ &\quad\times\left\{\varepsilon^{1/(4r)}(t-s)^{1/(2r)}+\varepsilon^{1/2}\left(\left\lfloor\frac{t}{\varepsilon}\right\rfloor-\left\lfloor\frac{s}{\varepsilon}\right\rfloor\right)^{(r-1)/(2r)} + \varepsilon^{(r+1)/(4r)}\left(\left\lfloor\frac{t}{\varepsilon}\right\rfloor-\left\lfloor\frac{s}{\varepsilon}\right\rfloor\right)^{1/4}\right\}.\end{aligned}\tag{35}$$

Moreover, we find from (35) that

$$\sup_{0\le\tau\le T}\hat{E}^{1/(2r)}\{\|\hat{u}_\tau^n\|_{L^2(\gamma)}^{2r}\}\le\frac{c(r,\Xi,d_1,\alpha,\|\Gamma\|,\|h'h\|_\infty,T)}{\varepsilon^{(r-1)/(4r)}n^{1/2}}.\tag{36}$$

(b) *Estimation of* $A_2$.

In the sequel, we use $*$ to denote the convolution of functions. By our assumption that $h\in\mathcal{S}(\mathbb{R}^{d_1})$, one finds that $\varrho_k^\varepsilon\in\mathcal{S}(\mathbb{R}^{d_1})$ for $k\in\mathbb{N}$. We define the function

$$\psi_\tau(\theta)\doteq\exp\{\ell(\theta)\tau\}\qquad\forall\,\theta\in\mathbb{R}^{d_1},\tau\in\mathbb{R}$$

and the operators

$$A_kf\doteq(1+\mathcal{B}_k)f,\qquad\mathcal{B}_kf\doteq\hat{\varrho}_k^\varepsilon*f,$$

$$
\begin{aligned}
&\hat{\varrho}^{\varepsilon}_k(\theta) \doteq \int_{\mathbb{R}^{d_1}} e_{-\theta}(x)\varrho^{\varepsilon}_k(x)\,dx, \qquad k\in\mathbb{N},\\
&T^{\lfloor t/\varepsilon\rfloor}_{t,s} f \doteq \psi_{t-\lfloor t/\varepsilon\rfloor\varepsilon} A_{\lfloor t/\varepsilon\rfloor}\left(\prod_{k=\lceil s/\varepsilon\rceil+1}^{\lfloor t/\varepsilon\rfloor-1}\{\psi_\varepsilon A_k\}\right)\psi_\varepsilon \mathcal{B}_{\lceil s/\varepsilon\rceil} f,\\
&T_{t,s} f \doteq \psi_{t-\lfloor t/\varepsilon\rfloor\varepsilon} A_{\lfloor t/\varepsilon\rfloor}\left(\prod_{k=\lceil s/\varepsilon\rceil+1}^{\lfloor t/\varepsilon\rfloor-1}\{\psi_\varepsilon A_k\}\right)\psi_{(\lceil s/\varepsilon\rceil+1)\varepsilon-s} f
\end{aligned}
\tag{37}
$$

for $f\in L^2(\mathbb{R}^{d_1};\gamma(d\theta))$, with the interpretations that the products go from right to left as one goes from the bottom. We find by (22) and (24)–(27) that

$$
\hat{v}^n_t(\theta) = \int_0^t (\ell(\theta) + (D\hat{Y}^{\varepsilon}_s * \hat{v}^n_{s-})(\theta))\,ds + \int_0^t (D\hat{Y}^{\varepsilon}_s * \hat{u}^n_{s-})(\theta)\,ds.
$$

Hence, $\hat{v}^n_t(\theta)$ is given by

$$
\hat{v}^n_t(\theta) = \sum_{k=1}^{\lfloor t/\varepsilon\rfloor} T^{\lfloor t/\varepsilon\rfloor}_{t,k\varepsilon}\hat{u}^n_{k\varepsilon-}(\theta). \tag{38}
$$

Moreover, we find for $0\le s<t\le T$ that

$$
\begin{aligned}
|\hat{v}^n_t(\theta)-\hat{v}^n_s(\theta)| &\le |(T_{t,s}-\psi_{t-s})\hat{v}^n_s(\theta)| + |(\psi_{t-s}-1)\hat{v}^n_s(\theta)|\\
&\quad + \left|\sum_{k=\lfloor s/\varepsilon\rfloor+1}^{\lfloor t/\varepsilon\rfloor} T^{\lfloor t/\varepsilon\rfloor}_{t,k\varepsilon}\hat{u}^n_{k\varepsilon-}(\theta)\right|.
\end{aligned}
\tag{39}
$$

Yet, recalling (37) and defining

$$
\tilde{T}^l_{t,s} f(\theta) = \psi_{t-l\varepsilon}\mathcal{B}_l\psi_\varepsilon T_{(l-1)\varepsilon,s} f(\theta) \qquad \forall\, l = \left\lceil\frac{s}{\varepsilon}\right\rceil+1,\ldots,\left\lfloor\frac{t}{\varepsilon}\right\rfloor, \tag{40}
$$

we see that, for any $\theta\in\mathbb{R}^{d_1}$, $0\le s\le T$,

$$
(T_{t,s}-\psi_{t-s})\hat{v}^n_s(\theta) = \sum_{l=\lceil s/\varepsilon\rceil+1}^{\lfloor t/\varepsilon\rfloor}\tilde{T}^l_{t,s}\hat{v}^n_s(\theta)
$$

and for any $\theta\in\mathbb{R}^{d_1}$, $0\le u<v\le T$,

$$
\sum_{k=\lfloor u/\varepsilon\rfloor+1}^{\lfloor v/\varepsilon\rfloor} T^{\lfloor v/\varepsilon\rfloor}_{v,k\varepsilon}\hat{u}^n_{k\varepsilon-}(\theta)
$$

are sums of, respectively, forward martingale and backward martingale differences. Thus, we find that

$$\hat{E}\{|(T_{t,s}-\psi_{t-s})\hat{v}_s^n(\theta)|^2\} = \sum_{l=\lfloor s/\varepsilon\rfloor+1}^{\lfloor t/\varepsilon\rfloor} \hat{E}\{|\tilde{T}_{t,s}^l \hat{v}_s^n(\theta)|^2\}, \tag{41}$$

$$\hat{E}\left\{\left|\sum_{k=\lfloor u/\varepsilon\rfloor+1}^{\lfloor v/\varepsilon\rfloor} T_{v,k\varepsilon}^{\lfloor v/\varepsilon\rfloor}\hat{u}_{k\varepsilon-}^n(\theta)\right|^2\right\} = \sum_{k=\lfloor u/\varepsilon\rfloor+1}^{\lfloor v/\varepsilon\rfloor} \hat{E}\{|T_{v,k\varepsilon}^{\lfloor v/\varepsilon\rfloor}\hat{u}_{k\varepsilon-}^n(\theta)|^2\}. \tag{42}$$

For $\varpi \in \mathbb{R}^{d_1}$, we define

$$m_\varpi(\cdot) \doteq (1+|\cdot|^2)^{-\gamma/2}(1+|\cdot+\varpi|^2)^{\gamma/2}.$$

Then, by Minkowski's integral inequality, classical multiplier theorem (see [18], page 96, Theorem 3), Jensen's inequality, independence, the assumption on $h$, (11) and Lemma 3, we find that

$$\begin{aligned}
\hat{E}\{\|\mathcal{B}_l f\|_{L^2(\gamma)}^2\} &= \hat{E}\left\{\int_{\mathbb{R}^{d_1}}\left|\int_{\mathbb{R}^{d_1}} \hat{\varrho}_l^\varepsilon(\varpi) f(\theta-\varpi)(1+|\theta|^2)^{\gamma/2}\,d\varpi\right|^2 d\theta\right\}\\
&\le \hat{E}\left\{\left(\int_{\mathbb{R}^{d_1}} |\hat{\varrho}_l^\varepsilon(\varpi)|\cdot\|f\cdot(1+|\cdot+\varpi|^2)^{\gamma/2}\|_{L^2}\,d\varpi\right)^2\right\}\\
&= \hat{E}\left\{\left(\int_{\mathbb{R}^{d_1}} |\hat{\varrho}_l^\varepsilon(\varpi)|\cdot\|m_\varpi\cdot f\cdot(1+|\cdot|^2)^{\gamma/2}\|_{L^2}\,d\varpi\right)^2\right\}\\
&\le c(d_1)\hat{E}\left\{\left(\int_{\mathbb{R}^{d_1}} |\hat{\varrho}_l^\varepsilon(\varpi)|(1+|\varpi|^2)^{[d_1/2]+1-\gamma/2}\,d\varpi\right)^2\|f\|_{L^2(\gamma)}^2\right\} \qquad (43)\\
&\le c(d_1)\hat{E}\{\|(1+|\cdot|^2)^{([d_1-2\gamma]+2)/2}\hat{\varrho}_l^\varepsilon\|_{L^2}^2\cdot\|f\|_{L^2(\gamma)}^2\}\\
&\le c(d_1)\hat{E}\{|\!|\!|\varrho_l^\varepsilon|\!|\!|_{[d_1-2\gamma]+2}^2\}\hat{E}\{\|f\|_{L^2(\gamma)}^2\}\\
&\le c(d_1,\langle\!\langle h\rangle\!\rangle_{[d_1-2\gamma]+2})\varepsilon\hat{E}\{\|f\|_{L^2(\gamma)}^2\}
\end{aligned}$$

for any $f\in L^2(\hat{\Omega},\ \mathcal{F}^{t_l-},\ L^2(\mathbb{R}^{d_1};\ \gamma(d\theta)))$, where

$$|\!|\!|\varrho_l^\varepsilon|\!|\!|_{[d_1-2\gamma]+2} = \left(\sum_{|\tau|\le[d_1-2\gamma]+2}\|D^\tau\varrho_l^\varepsilon\|_{L^2}^2\right)^{1/2}$$

is the standard Sobolev $W^{[d_1-2\gamma]+2,2}$-norm of $\varrho_l^\varepsilon$. Moreover, using (40), the fact that $|\psi_\varepsilon(\theta)|\le 1$, independence, (43) and recursion, we find that

$$\begin{aligned}
&\max_{\lfloor s/\varepsilon\rfloor+1\le l\le\lfloor t/\varepsilon\rfloor} \hat{E}^{1/2}\{\|\tilde{T}_{t,s}^l\hat{v}_s^n\|_{L^2(\gamma)}^2\}\\
&\quad\le \max_{\lfloor s/\varepsilon\rfloor+1\le l\le\lfloor t/\varepsilon\rfloor}(c(d_1,\langle\!\langle h\rangle\!\rangle_{[d_1-2\gamma]+2})\varepsilon\hat{E}\{\|T_{(l-1)\varepsilon,s}\hat{v}_s^n\|_{L^2(\gamma)}^2\})^{1/2}
\end{aligned}$$

$$
\begin{aligned}
&= \max_{\lfloor s/\varepsilon\rfloor+1\le l\le\lfloor t/\varepsilon\rfloor}\biggl(c(d_1,\langle\!\langle h\rangle\!\rangle_{[d_1-2\gamma]+2})\varepsilon \\
&\qquad\times\int_{\mathbb{R}^{d_1}}\hat{E}\{|\psi_\varepsilon T_{(l-2)\varepsilon,s}\hat{v}_s^n(\theta)|^2 \\
&\qquad\qquad+|\mathcal{B}_{l-1}\psi_\varepsilon T_{(l-2)\varepsilon,s}\hat{v}_s^n(\theta)|^2\}\gamma(d\theta)\biggr)^{1/2} \\
&\le \max_{\lfloor s/\varepsilon\rfloor+1\le l\le\lfloor t/\varepsilon\rfloor}(c(d_1,\langle\!\langle h\rangle\!\rangle_{[d_1-2\gamma]+2})\varepsilon(1+c(d_1,\langle\!\langle h\rangle\!\rangle_{[d_1-2\gamma]+2})\varepsilon) \\
&\qquad\times\hat{E}\{\|T_{(l-2)\varepsilon,s}\hat{v}_s^n\|^2_{L^2(\gamma)}\})^{1/2} \\
&\le \max_{\lfloor s/\varepsilon\rfloor+1\le l\le\lfloor t/\varepsilon\rfloor}(c(d_1,\langle\!\langle h\rangle\!\rangle_{[d_1-2\gamma]+2})\varepsilon(1+c(d_1,\langle\!\langle h\rangle\!\rangle_{[d_1-2\gamma]+2})\varepsilon)^{l-1} \\
&\qquad\times\hat{E}\{\|\hat{v}_s^n\|^2_{L^2(\gamma)}\})^{1/2} \qquad (44) \\
&\le c(d_1,\langle\!\langle h\rangle\!\rangle_{[d_1-2\gamma]+2})\varepsilon^{1/2}\hat{E}^{1/2}\{\|\hat{v}_s^n\|^2_{L^2(\gamma)}\}.
\end{aligned}
$$

Now, in a similar manner to (44), we find from (36) with $r=2$ and Jensen's inequality that

$$
\begin{aligned}
&\max_{\lfloor u/\varepsilon\rfloor+1\le k\le\lfloor v/\varepsilon\rfloor}\hat{E}^{1/2}\{\|T^{\lfloor v/\varepsilon\rfloor}_{v,k\varepsilon}\hat{u}^n_{k\varepsilon-}\|^2_{L^2(\gamma)}\} \\
&\qquad\le c(d_1,\langle\!\langle h\rangle\!\rangle_{[d_1-2\gamma]+2})\varepsilon^{1/2}\hat{E}^{1/2}\{\|\hat{u}_s^n\|^2_{L^2(\gamma)}\} \qquad (45) \\
&\qquad\le \frac{c(\Xi,d_1,\alpha,\|\Gamma\|,\langle\!\langle h\rangle\!\rangle_{[d_1-2\gamma]+2},T)\varepsilon^{3/8}}{n^{1/2}}.
\end{aligned}
$$

Hence, combining (42), (45) and (38), we find that

$$
\begin{aligned}
&\hat{E}^{1/2}\Biggl\{\Biggl\|\sum_{k=\lfloor s/\varepsilon\rfloor+1}^{\lfloor t/\varepsilon\rfloor}T^{\lfloor t/\varepsilon\rfloor}_{t,k\varepsilon}\hat{u}^n_{k\varepsilon-}\Biggr\|^2_{L^2(\gamma)}\Biggr\} \\
&\qquad\le\frac{c(\Xi,d_1,\alpha,\|\Gamma\|,\langle\!\langle h\rangle\!\rangle_{[d_1-2\gamma]+2},T)\varepsilon^{3/8}}{n^{1/2}}\biggl(\biggl\lfloor\frac{t}{\varepsilon}\biggr\rfloor-\biggl\lfloor\frac{s}{\varepsilon}\biggr\rfloor\biggr)^{1/2}, \qquad (46)
\end{aligned}
$$

$$
\hat{E}^{1/2}\{\|\hat{v}_t^n\|^2_{L^2(\gamma)}\}\le\frac{c(\Xi,d_1,\alpha,\|\Gamma\|,\langle\!\langle h\rangle\!\rangle_{[d_1-2\gamma]+2},T)}{\varepsilon^{1/8}n^{1/2}}. \qquad (47)
$$

Replacing $\gamma(d\theta)$ with $(|\theta|^\alpha\ln|\theta|)^2\gamma(d\theta)$, noting that $\gamma<-(d_1/2+2\alpha)$ by the assumption and repeating the above arguments, one finds that

$$
\begin{aligned}
&\hat{E}^{1/2}\biggl\{\int_{\mathbb{R}^d}(|\theta|^\alpha\ln|\theta|)^2|\hat{v}_s^n|^2(\theta)\gamma(d\theta)\biggr\} \qquad (48) \\
&\qquad\le\frac{c(\Xi,d_1,\alpha,\|\Gamma\|,\langle\!\langle h\rangle\!\rangle_{[d_1-2\gamma]+2},T)}{\varepsilon^{1/8}n^{1/2}}.
\end{aligned}
$$

Now, it follows by (41), (44) and (**??**) that

$$\begin{aligned}\hat{E}^{1/2}\{\|(T_{t,s}-\psi_{t-s})\hat{v}_s\|^2_{L^2(\gamma)}\}\\ \leq \frac{c(\Xi,d_1,\alpha,\|\Gamma\|,\langle\langle h\rangle\rangle_{[d_1-2\gamma]+2},T)\varepsilon^{3/8}}{n^{1/2}}\left(\left\lfloor\frac{t}{\varepsilon}\right\rfloor-\left\lfloor\frac{s}{\varepsilon}\right\rfloor\right)^{1/2}.\end{aligned} \tag{49}$$

Finally, using the bound $|\psi_{t-s}(\theta)-1|^2 \leq c(\alpha,\|\Gamma\|)(|\theta|^\alpha \ln|\theta|)^2|t-s|^2$, (39), (46), (48) and (49), one finds that

$$\begin{aligned}&\hat{E}^{1/2}\{\|\hat{v}^n_t-\hat{v}^n_s\|^2_{L^2(\gamma)}\}\\ &\quad\leq \frac{c(\Xi,d_1,\alpha,\|\Gamma\|,\langle\langle h\rangle\rangle_{[d_1-2\gamma]+2},T)}{\varepsilon^{1/8}n^{1/2}}\\ &\qquad\times\left\{(t-s)+\varepsilon^{1/2}\left(\left\lfloor\frac{t}{\varepsilon}\right\rfloor-\left\lfloor\frac{s}{\varepsilon}\right\rfloor\right)^{1/2}\right\}.\end{aligned} \tag{50}$$

(c) *Estimation of* $A_3$.

Note that the solution $\hat{\chi}^n_t(\theta)$ defined by (26) and (25) can be written as

$$\hat{\chi}^n_t(\theta)=\psi_{t-\lfloor t/\varepsilon\rfloor\varepsilon}\prod_{k=1}^{\lfloor t/\varepsilon\rfloor}\{A_k\psi_\varepsilon\}\langle\mu^n_0-\mu_0,e_{-\theta}\rangle$$

and

$$\hat{E}^{1/2}\{|\langle\mu^n_0-\mu_0,e_{-\theta}\rangle|^2\}\leq\frac{4}{n^{1/2}}\qquad\forall\, n\in\mathbb{N},\theta\in\mathbb{R}^{d_1}. \tag{51}$$

Then, one finds similarly to (50) that

$$\begin{aligned}&\hat{E}^{1/2}\{\|\hat{\chi}^n_t-\hat{\chi}^n_s\|^2_{L^2(\gamma)}\}\\ &\quad\leq \frac{c(\Xi,d_1,\alpha,\|\Gamma\|,\langle\langle h\rangle\rangle_{[d_1-2\gamma]+2},T)\varepsilon^{1/2}}{n^{1/2}}\\ &\qquad\times\left\{(t-s)+\left(\left\lfloor\frac{t}{\varepsilon}\right\rfloor-\left\lfloor\frac{s}{\varepsilon}\right\rfloor\right)^{1/2}\right\}.\end{aligned} \tag{52}$$

Therefore, (12) follows from (28), (35) with $r=2$, Jensen's inequality, (50) and (52). By virtue of (51), (13) is an immediate consequence of (12) by letting $\Xi=\varepsilon^{1/2}$. □

3.3. *Proof of Corollary* 5. Since $\alpha=2$, we find by (30) and (18) that

$$\begin{aligned}&\hat{E}^{1/r}\{|\hat{u}^{n,e}_t(\theta)-\hat{u}^{n,e}_s(\theta)|^{2r}\}\\ &\quad\leq\frac{c(r)}{n}\sup_{0\leq\tau\leq T}\hat{E}^{1/r}\{\langle\mu^n_\tau,1\rangle^r\}\end{aligned}$$

$$
\begin{aligned}
&\times\biggl\{\biggl(\exp\biggl\{-(t-s)\int_{S_{d_1}}|\theta' z|^2\Gamma(dz)\biggr\}-1\biggr)^{2r}\\
&\qquad\times\biggl(\int_0^s\exp\biggl\{-2(s-\tau)\int_{S_{d_1}}|\theta' z|^2\Gamma(dz)\biggr\}d\biggl(\tau\int_{S_{d_1}}|\theta' z|^2\Gamma(dz)\biggr)\biggr)^r\\
&\qquad+\biggl(\int_s^t\exp\biggl\{-2(t-\tau)\int_{S_{d_1}}|\theta' z|^2\Gamma(dz)\biggr\}d\biggl(\tau\int_{S_{d_1}}|\theta' z|^2\Gamma(dz)\biggr)\biggr)^r\biggr\}^{1/r}\\
&\le\frac{c(r)}{n}\sup_{0\le\tau\le T}\hat{E}^{1/r}\{\langle\mu_\tau^n,1\rangle^r\}\biggl(1-\exp\biggl\{-2(t-s)\int_{S_{d_1}}|\theta' z|^2\Gamma(dz)\biggr\}\biggr)\\
&\le\frac{c(r,\|\Gamma\|)}{n}\sup_{0\le\tau\le T}\hat{E}^{1/r}\{\langle\mu_\tau^n,1\rangle^r\}(t-s)|\theta|^2.
\end{aligned}
\tag{53}
$$

Replacing (31) with (53), we find similarly to (35) that

$$
\begin{aligned}
\hat{E}^{1/(2r)}\{\|\hat{u}_t^n-\hat{u}_s^n\|_{L^2(\gamma)}^{2r}\}&\le\frac{c(r,\Xi,d_1,\|\Gamma\|,\|h'h\|_\infty,T)}{\varepsilon^{1/4}n^{1/2}}\\
&\quad\times\biggl\{\varepsilon^{1/(4r)}(t-s)^{1/2}+\varepsilon^{1/2}\biggl(\biggl\lfloor\frac{t}{\varepsilon}\biggr\rfloor-\biggl\lfloor\frac{s}{\varepsilon}\biggr\rfloor\biggr)^{(r-1)/(2r)}\\
&\qquad+\varepsilon^{(r+1)/(4r)}\biggl(\biggl\lfloor\frac{t}{\varepsilon}\biggr\rfloor-\biggl\lfloor\frac{s}{\varepsilon}\biggr\rfloor\biggr)^{1/4}\biggr\}.
\end{aligned}
$$

Letting $r=2$, we then find that, for $\beta>1/8$,

$$
\begin{aligned}
\hat{E}\{\|\hat{u}_t^n-\hat{u}_s^n\|_{L^2(\gamma)}^4\}&\le\frac{c(\Xi,d_1,\|\Gamma\|,\|h'h\|_\infty,T)}{\varepsilon n^2}\\
&\quad\times\biggl\{\varepsilon^{1/2}(t-s)^2+\varepsilon^2\biggl(\biggl\lfloor\frac{t}{\varepsilon}\biggr\rfloor-\biggl\lfloor\frac{s}{\varepsilon}\biggr\rfloor\biggr)+\varepsilon^{3/2}\biggl(\biggl\lfloor\frac{t}{\varepsilon}\biggr\rfloor-\biggl\lfloor\frac{s}{\varepsilon}\biggr\rfloor\biggr)\biggr\}\\
&\le\frac{c(\Xi,d_1,\|\Gamma\|,\|h'h\|_\infty,T)}{\varepsilon n^2}\\
&\quad\times\biggl\{\varepsilon^{1/2}(t-s)^2+\varepsilon^2\biggl(\biggl\lfloor\frac{t}{\varepsilon}\biggr\rfloor-\biggl\lfloor\frac{s}{\varepsilon}\biggr\rfloor\biggr)^{(1/2)+4\beta}\\
&\qquad+\varepsilon^{3/2}\biggl(\biggl\lfloor\frac{t}{\varepsilon}\biggr\rfloor-\biggl\lfloor\frac{s}{\varepsilon}\biggr\rfloor\biggr)^{(1/2)+4\beta}\biggr\}.
\end{aligned}
$$

Thus, by Lemma 10, we find that

$$
\begin{aligned}
\hat{E}^{1/2}\biggl\{\sup_{0\le s<t\le T}\|\hat{u}_t^n-\hat{u}_s^n\|_{L^2(\gamma)}^2\biggr\}&\le\hat{E}^{1/4}\biggl\{\sup_{0\le s<t\le T}\|\hat{u}_t^n-\hat{u}_s^n\|_{L^2(\gamma)}^4\biggr\}\\
&\le\frac{c(\Xi,d_1,\beta,\|\Gamma\|,\|h'h\|_\infty,T)}{\varepsilon^\beta n^{1/2}}.
\end{aligned}
\tag{54}
$$

Similarly, by (50), (52) and Lemma 10, we find that

$$\hat{E}^{1/2}\Big\{\sup_{0\le s<t\le T}\|\hat{v}_t^n-\hat{v}_s^n\|^2_{L^2(\gamma)}\Big\}\le\frac{c(\Xi,d_1,\beta,\|\Gamma\|,\langle\!\langle h\rangle\!\rangle_{[d_1-2\gamma]+2},T)}{\varepsilon^\beta n^{1/2}} \tag{55}$$

and

$$\hat{E}^{1/2}\Big\{\sup_{0\le s<t\le T}\|\hat{\chi}_t^n-\hat{\chi}_s^n\|^2_{L^2(\gamma)}\Big\}\le\frac{c(\Xi,d_1,\beta,\|\Gamma\|,\langle\!\langle h\rangle\!\rangle_{[d_1-2\gamma]+2},T)}{\varepsilon^\beta n^{1/2}}. \tag{56}$$

Therefore, (14) follows from (28), (54), (55) and (56). By virtue of (51), (15) is an immediate consequence of (14) by letting $\Xi=\varepsilon^{1/2}$. □

## APPENDIX: PROOFS OF LEMMA 1, PROPOSITION 2 AND LEMMA 3

In the current section we give the proofs of Lemma 1, Proposition 2 and Lemma 3. We realize that similar results are well known in a variety of settings and only give them for the sake of completeness.

Proof of Lemma 1. For $\varphi\in\mathcal{T}$, we have that

$$\varphi(X_t)-\varphi(X_0)=\int_0^t\mathcal{L}\varphi(X_s)\,ds+\mathcal{M}_t(\varphi),$$

where $\mathcal{M}_t(\varphi)$ is an $\mathcal{X}^t$-martingale. Then,

$$\begin{aligned}\varphi(X_t)\eta_t^{-1}&=\varphi(X_0)+\int_0^t\mathcal{L}\varphi(X_s)\eta_s^{-1}\,ds+\int_0^t\eta_{s-}^{-1}\,d\mathcal{M}_s(\varphi)\\&\quad+\sum_{k=1}^{\lfloor t/\varepsilon\rfloor}\varphi(X_{t_k})\eta_{t_k}^{-1}\bigg(1-\exp\bigg\{-(Y_{t_k}-Y_{t_{k-1}})'h(X_{t_k})\\&\qquad+\frac{(t_k-t_{k-1})(h'h)(X_{t_k})}{2}\bigg\}\bigg).\end{aligned}$$

By the independence of $X$ and $Y$ under $P$, we find that $\mathcal{M}_t(\varphi)$ is also an $\mathcal{X}^t\vee\mathcal{Y}^T$-martingale so $E\{\int_0^t\eta_{s-}^{-1}\,d\mathcal{M}_s(\varphi)|\mathcal{Y}^t\}=0$. Hence,

$$\begin{aligned}\langle\mu_t,\varphi\rangle&=\langle\mu_0,\varphi\rangle+\int_0^t\langle\mu_s,\mathcal{L}\varphi\rangle\,ds\\&\quad+\sum_{k=1}^{\lfloor t/\varepsilon\rfloor}\bigg\langle\mu_{t_k},\varphi\bigg(1-\exp\bigg\{-(Y_{t_k}-Y_{t_{k-1}})'h+\frac{(t_k-t_{k-1})h'h}{2}\bigg\}\bigg)\bigg\rangle.\end{aligned} \tag{57}$$

On the other hand, we obtain from the definition of $\mu_t$ and the stochastic continuity of $X$ that, for any continuous bounded function $\varphi$ on $\mathbb{R}^{d_1}$, and $k \geq 1$,

$$\left\langle \mu_{t_k}, \varphi\left(1-\exp\left\{-(Y_{t_k}-Y_{t_{k-1}})'h+\frac{(t_k-t_{k-1})h'h}{2}\right\}\right)\right\rangle$$

$$\begin{aligned}
&= E^X\left\{\varphi(X_{t_k})\left(1-\exp\left\{-(Y_{t_k}-Y_{t_{k-1}})'h(X_{t_k})+\frac{(t_k-t_{k-1})(h'h)(X_{t_k})}{2}\right\}\right)\right.\\
&\qquad\qquad \left.\times\prod_{l=1}^{k}\exp\left\{(Y_{t_l}-Y_{t_{l-1}})'h(X_{t_l})-\frac{(t_l-t_{l-1})(h'h)(X_{t_l})}{2}\right\}\right\}\\
&= E^X\left\{\varphi(X_{t_k})\left(\exp\left\{(Y_{t_k}-Y_{t_{k-1}})'h(X_{t_k})-\frac{(t_k-t_{k-1})(h'h)(X_{t_k})}{2}\right\}-1\right)\right.\\
&\qquad\qquad \left.\times\prod_{l=1}^{k-1}\exp\left\{(Y_{t_l}-Y_{t_{l-1}})'h(X_{t_l})-\frac{(t_l-t_{l-1})(h'h)(X_{t_l})}{2}\right\}\right\} \qquad (58)\\
&= \lim_{t\uparrow t_k}\left\langle \mu_t, \varphi\left(\exp\left\{(Y_{t_k}-Y_{t_{k-1}})'h-\frac{(t_k-t_{k-1})h'h}{2}\right\}-1\right)\right\rangle\\
&= \left\langle \mu_{t_k-}, \varphi\left(\exp\left\{(Y_{t_k}-Y_{t_{k-1}})'h-\frac{(t_k-t_{k-1})h'h}{2}\right\}-1\right)\right\rangle,
\end{aligned}$$

where $E^X$ is the expectation taken only with respect to $X$. Further, (58) holds for any $\varphi \in \mathcal{B}_b(\mathbb{R}^{d_1})$ by the monotone class theorem. Substituting (58) into (57), we get (4).

The uniqueness of $\mu_t$ can be proved by the action of $\mathcal{L}$ on the trigonometric polynomials and induction. In fact, suppose that $\{\mu_t, t \geq 0\}$ and $\{\nu_t, t \geq 0\}$ satisfy (4), and $\mu_t = \nu_t$ for $t \leq t_k$ for some $k \geq 0$. Note that

$$\mathcal{L}e_{-\theta} = \begin{cases} -\left(\displaystyle\int_{S_{d_1}} |\theta' z|^\alpha\left(1+i\,\mathrm{sign}(\theta' z)\tan\left(\frac{\alpha\pi}{2}\right)\right)\Gamma(dz)\right)e_{-\theta}, & \text{for } \alpha \neq 1,\\ -\left(\displaystyle\int_{S_{d_1}} |\theta' z|\left(1-\frac{2i}{\pi}\mathrm{sign}(\theta' z)\ln|\theta' z|\right)\Gamma(dz)\right)e_{-\theta}, & \text{for } \alpha = 1.\end{cases}$$

From (4), one finds that, for any $\theta \in \mathbb{R}^{d_1}$, $t_k \leq t < t_{k+1}$,

$$\begin{aligned}
&\langle \mu_t, e_{-\theta}\rangle\\
&\quad= \begin{cases} \langle \mu_{t_k}, e_{-\theta}\rangle \exp\left\{-(t-t_k)\displaystyle\int_{S_{d_1}} |\theta' z|^\alpha\left(1+i\,\mathrm{sign}(\theta' z)\tan\left(\frac{\alpha\pi}{2}\right)\right)\Gamma(dz)\right\}, \\ \qquad\qquad \text{for } \alpha \neq 1,\\ \langle \mu_{t_k}, e_{-\theta}\rangle \exp\left\{-(t-t_k)\displaystyle\int_{S_{d_1}} |\theta' z|\left(1-\frac{2i}{\pi}\mathrm{sign}(\theta' z)\ln|\theta' z|\right)\Gamma(dz)\right\}, \\ \qquad\qquad \text{for } \alpha = 1\end{cases}\\
&\quad= \langle \nu_t, e_{-\theta}\rangle.
\end{aligned}$$

Since the set of trigonometric polynomials is measure-determining, $\mu_t = \nu_t$, $t_k \le t < t_{k+1}$. Hence, by (4), we find that $\mu_{t_{k+1}} = \nu_{t_{k+1}}$. Therefore, the uniqueness of $\mu_t$ follows by induction. □

PROOF OF LEMMA 3. Let $W$ be a standard $\mathbb{R}^{d_2}$-valued Brownian motion on $(\hat{\Omega}, \hat{\mathcal{F}}, \hat{P})$. We fix an $x \in \mathbb{R}^{d_1}$ and define

$$Z_t^x \doteq \exp\left\{(W_t)'h(x) - \frac{t(h'h)(x)}{2}\right\}.$$

Then,

$$\hat{E}\{|\varrho_k^\varepsilon(x)|^r\} = \hat{E}\{|Z_\varepsilon^x - 1|^r\}$$

for any $k \in \mathbb{N}$. By Burkholder's, Minkowski's integral and Jensen's inequalities, we find that, for any $0 \le t \le \varepsilon$,

$$\begin{aligned}
\hat{E}\{|Z_t^x - 1|^r\} &\le c(r)\|h'h\|_\infty^{r/2}\hat{E}\left\{\left(\int_0^t (Z_s^x)^2\,ds\right)^{r/2}\right\} \\
&\le c(r)\|h'h\|_\infty^{r/2}\left(\hat{E}\left\{\left(\int_0^t (Z_s^x - 1)^2\,ds\right)^{r/2}\right\} + \varepsilon^{r/2}\right) \\
&\le c(r)\|h'h\|_\infty^{r/2}\left(\left(\int_0^t (\hat{E}\{|Z_s^x - 1|^r\})^{2/r}\,ds\right)^{r/2} + \varepsilon^{r/2}\right) \\
&\le c(r)\|h'h\|_\infty^{r/2}\left(\varepsilon^{(r/2)-1}\int_0^t \hat{E}\{|Z_s^x - 1|^r\}\,ds + \varepsilon^{r/2}\right),
\end{aligned}$$

where we have assumed without loss of generality that $r \ge 2$ above. Applying Gronwall's inequality, one then discovers that

$$\sup_{0 \le t \le \varepsilon} \hat{E}\{|Z_t^x - 1|^r\} \le c(r)\|h'h\|_\infty^{r/2}\varepsilon^{r/2}$$

and the lemma follows. □

PROOF OF PROPOSITION 2. To ease the notation in the sequel, we let $\xi_k = \xi_k^\varepsilon$. For $\varphi \in \mathcal{T}$, we define

$$\begin{aligned}
\mathcal{M}_t^n(\varphi) \doteq{}& \frac{1}{n}\sum_{k=0}^{\lfloor t/\varepsilon \rfloor}\sum_{i=1}^{\|\mu_{t_k}^n\|}\left(\varphi(X_{t_{k+1}\wedge t}^{i,n}) - \varphi(X_{t_k}^{i,n}) - \int_{t_k}^{t_{k+1}\wedge t}\mathcal{L}\varphi(X_u^{i,n})\,du\right) \\
&+ \frac{1}{n}\sum_{k=1}^{\lfloor t/\varepsilon \rfloor}\sum_{i=1}^{\|\mu_{t_k-}^n\|}\langle \delta_{X_{t_k-}^{i,n}}, \varphi\rangle \qquad (59) \\
&\times (\operatorname{sign}(\xi_k(X_{t_k-}^{i,n}))\mathbb{1}_{\{U^{i,k} \in [0, |\xi_k(X_{t_k-}^{i,n})|)\}} - \xi_k(X_{t_k-}^{i,n})).
\end{aligned}$$

Then, we find from our algorithm and (1) that (6) holds. Recalling that the $\{U^{i,k}\}$ are independent and compensating the square of the jumps in the second term of (59), we find that $\{\mathcal{M}_t^n(\varphi)\}_{t\geq 0}$ is a cádlág $\{\mathcal{G}^t\}_{t\geq 0}$-martingale satisfying (7).

Now, turning to bounding the difference between the quadratic variation $[\mathcal{M}^n(\varphi)]_t$ and the expected quadratic variation $\hat{E}^U\{[\mathcal{M}^n(\varphi)]_t\}$, we define

$$A^{i,k} \doteq \frac{1}{n}\langle \delta_{X^{i,n}_{t_k-}}, \varphi\rangle (\operatorname{sign}(\xi_k(X^{i,n}_{t_k-}))1_{\{U^{i,k}\in[0,|\xi_k(X^{i,n}_{t_k-})|)\}} - \xi_k(X^{i,n}_{t_k-})).$$

Letting $\{f_k\}_{k=1}^\infty \subset \mathbb{R}$, recognizing the martingale transform and using Burkholder's and Jensen's inequalities, we bound

$$\begin{aligned}
&\hat{E}\left\{\left|\sum_{k=\lfloor s/\varepsilon\rfloor+1}^{\lfloor t/\varepsilon\rfloor} f_k([\mathcal{M}^n(\varphi)]_{k\varepsilon} - \hat{E}^U\{[\mathcal{M}^n(\varphi)]_{k\varepsilon}\})\right|^r\right\}\\
&\quad = \hat{E}\left\{\left|\sum_{k=\lfloor s/\varepsilon\rfloor+1}^{\lfloor t/\varepsilon\rfloor} f_k\left(\left(\sum_{i=1}^{\|\mu^n_{t_k-}\|} A^{i,k}\right)^2 - \hat{E}^U\left\{\left(\sum_{i=1}^{\|\mu^n_{t_k-}\|} A^{i,k}\right)^2\right\}\right)\right|^r\right\}\\
&\quad \le c(r)\hat{E}\left\{\left|\sum_{k=\lfloor s/\varepsilon\rfloor+1}^{\lfloor t/\varepsilon\rfloor} f_k^2\left(\left(\sum_{i=1}^{\|\mu^n_{t_k-}\|} A^{i,k}\right)^2 - \hat{E}^U\left\{\left(\sum_{i=1}^{\|\mu^n_{t_k-}\|} A^{i,k}\right)^2\right\}\right)^2\right|^{r/2}\right\} \qquad (60)\\
&\quad \le c(r)\left(\sum_{k=\lfloor s/\varepsilon\rfloor+1}^{\lfloor t/\varepsilon\rfloor} f_k^2\right)^{(r/2)-1}\\
&\qquad\times \sum_{k=\lfloor s/\varepsilon\rfloor+1}^{\lfloor t/\varepsilon\rfloor} f_k^2 \hat{E}\left\{\left|\left(\sum_{i=1}^{\|\mu^n_{t_k-}\|} A^{i,k}\right)^2 - \hat{E}^U\left\{\left(\sum_{i=1}^{\|\mu^n_{t_k-}\|} A^{i,k}\right)^2\right\}\right|^r\right\}.
\end{aligned}$$

However, defining the filtrations $\{\mathcal{F}^m_{k,+}\}_{m=1}^\infty$ and $\{\mathcal{F}^m_{k,-}\}_{m=1}^\infty$ via

$$\mathcal{F}^m_{k,+} \doteq \mathcal{G}^{t_k-} \vee \sigma\{U^{i,k}, i\le m\}, \qquad \mathcal{F}^m_{k,-} \doteq \mathcal{G}^{t_k-} \vee \sigma\{U^{i,k}, i\ge m\},$$

we find that

$$m \to \left(\sum_{i=1}^m A^{i,k}\right)^2 - \hat{E}^U\left\{\left(\sum_{i=1}^m A^{i,k}\right)^2\right\}$$

is an $\{\mathcal{F}^m_{k,+}\}_{m=1}^\infty$-martingale and $m \to A^{i,k}\sum_{j=m}^{i-1} A^{j,k}$ is a backward $\{\mathcal{F}^m_{k,-}\}_{m=1}^\infty$-martingale for each $i$. This means we can again apply Burkholder's, Jensen's and $2ab \le a^2+b^2$ inequalities and use the independence of $\{U^{i,k}\}$ to find that

$$\hat{E}\left\{\left|\left(\sum_{i=1}^{\|\mu^n_{t_k-}\|} A^{i,k}\right)^2 - \hat{E}^U\left\{\left(\sum_{i=1}^{\|\mu^n_{t_k-}\|} A^{i,k}\right)^2\right\}\right|^r \middle| \mathcal{G}^{t_k-}\right\}$$

$$\leq c(r)\|\mu^n_{t_k-}\|^{(r/2)-1}\sum_{i=1}^{\|\mu^n_{t_k-}\|}\hat{E}\Biggl\{|(A^{i,k})^2-\hat{E}^U\{(A^{i,k})^2\}|^r + \left|2A^{i,k}\sum_{j=1}^{i-1}A^{j,k}\right|^r\Bigg|\mathcal{G}^{t_k-}\Biggr\}$$

$$\leq c(r)\|\mu^n_{t_k-}\|^{r-2}\sum_{i=1}^{\|\mu^n_{t_k-}\|}\Biggl(\hat{E}\{|(A^{i,k})^2-\hat{E}^U\{(A^{i,k})^2\}|^r|\mathcal{G}^{t_k-}\} + \sum_{j=1}^{i-1}\hat{E}\{|A^{i,k}A^{j,k}|^r|\mathcal{G}^{t_k-}\}\Biggr) \tag{61}$$

$$\leq c(r)\|\mu^n_{t_k-}\|^{r-2}\sum_{i=1}^{\|\mu^n_{t_k-}\|}\Biggl(i\hat{E}^U\{|(A^{i,k})^2-\hat{E}^U\{(A^{i,k})^2\}|^r\} + (i-1)(\hat{E}^U\{(A^{i,k})^2\})^r + \sum_{j=1}^{i-1}(\hat{E}^U\{|(A^{j,k})^2-\hat{E}^U\{(A^{j,k})^2\}|^r\} + (\hat{E}^U\{(A^{j,k})^2\})^r)\Biggr)$$

$$\leq c(r)\|\mu^n_{t_k-}\|^{r-1}\sum_{i=1}^{\|\mu^n_{t_k-}\|}(\hat{E}^U\{|(A^{i,k})^2-\hat{E}^U\{(A^{i,k})^2\}|^r\} + (\hat{E}^U\{(A^{i,k})^2\})^r), \qquad \hat{P}\text{-a.s.}$$

Now, $\|\mu^n_{t_k-}\| = n\langle\mu^n_{t_k-},1\rangle$ and it follows by direct calculation of $\hat{E}^U\{|(A^{i,k})^2-\hat{E}^U\{(A^{i,k})^2\}|^r\}$ that

$$\hat{E}\Biggl\{(n\langle\mu^n_{t_k-},1\rangle)^{r-1}\sum_{i=1}^{\|\mu^n_{t_k-}\|}|(A^{i,k})^2-\hat{E}^U\{(A^{i,k})^2\}|^r\Biggr\}$$

$$= \hat{E}\biggl\{\frac{1}{n^r}\langle\mu^n_{t_k-},1\rangle^{r-1}\langle\mu^n_{t_k-},|\varphi|^{2r}\{|1-3|\xi_k|+2\xi_k^2|^r|\xi_k| + |2\xi_k^2-|\xi_k||^r(1-|\xi_k|)\}\rangle\biggr\}.$$

Next, conditioning on $\sigma\{\mu^n_{t_k-}\}$, using the independence of the increments of $Y$ and Lemma 3, we find that

$$\begin{aligned}
&\hat{E}\Biggl\{(n\langle\mu^n_{t_k-},1\rangle)^{r-1}\sum_{i=1}^{\|\mu^n_{t_k-}\|}|(A^{i,k})^2-\hat{E}^U\{(A^{i,k})^2\}|^r\Biggr\}\\
&\qquad\le\frac{c\|h'h\|_\infty^{1/2}\varepsilon^{1/2}}{n^r}\hat{E}\{\langle\mu^n_{t_k-},1\rangle^r\}\|\varphi\|_\infty^{2r}.
\end{aligned}\tag{62}$$

By Lemma 3 and the fact that

$$(\hat{E}^U\{(A^{i,k})^2\})^r=(|\xi_k|-\xi_k^2)^r n^{-2r}\varphi^{2r}(X^{i,n}_{t_k-}),$$

we find that

$$\begin{aligned}
&\hat{E}\Biggl\{(n\langle\mu^n_{t_k-},1\rangle)^{r-1}\sum_{i=1}^{\|\mu^n_{t_k-}\|}(\hat{E}^U\{(A^{i,k})^2\})^r\Biggr\}\\
&\qquad\le\frac{c(r)\|h'h\|_\infty^{r/2}\varepsilon^{r/2}}{n^r}\hat{E}\{\langle\mu^n_{t_k-},1\rangle^r\}\|\varphi\|_\infty^{2r}.
\end{aligned}\tag{63}$$

Then, substituting (62) and (63) into (61) and (60), we find that

$$\begin{aligned}
&\hat{E}\Biggl\{\Biggl|\sum_{k=\lfloor s/\varepsilon\rfloor+1}^{\lfloor t/\varepsilon\rfloor}f_k([\mathcal{M}^n(\varphi)]_{k\varepsilon}-\hat{E}^U\{[\mathcal{M}^n(\varphi)]_{k\varepsilon}\})\Biggr|^r\Biggr\}\\
&\quad\le c(r)(\|h'h\|_\infty\vee 1)^{r/2}\Biggl(\sum_{k=\lfloor s/\varepsilon\rfloor+1}^{\lfloor t/\varepsilon\rfloor}f_k^2\Biggr)^{r/2}\frac{\varepsilon^{1/2}}{n^r}\Bigl(\sup_{0\le\tau\le T}\hat{E}\{\langle\mu^n_\tau,1\rangle^r\}\Bigr)\|\varphi\|_\infty^{2r}
\end{aligned}$$

for some constant $c(r)>0$ independent of $d_1$, $d_2$, $\varepsilon$, $n$, $t$, $s$, $\varphi$. $\square$

**Acknowledgments.** The authors gratefully acknowledge the very helpful suggestions and comments of an anonymous referee, an Associate Editor and the Editor, which helped improve the presentation of this paper.

Department of Mathematical
and Statistical Sciences
University of Alberta
Edmonton, Alberta
Canada T6G 2G1
e-mail: mkouritz@math.ualberta.ca

Department of Mathematics
and Statistics
Concordia University
7141 Sherbrooke Street West
Montreal, Quebec
Canada H4B 1R6
e-mail: wsun@alcor.concordia.ca